\documentclass[dvips]{article}
\usepackage{graphicx}

\font\BBb=msbm10 at 12pt

\newcommand{\square}{\hbox{${\vcenter{\hrule height.4pt
  \hbox{\vrule width.4pt height6pt \kern6pt
     \vrule width.4pt}
  \hrule height.4pt}}$}}

\newcommand {\qed} {\hfill \nobreak \square \medbreak}

\newcommand{\be}{\begin{equation}}
      \newcommand{\ee}{\end{equation}}
      \newcommand{\ba}{\begin{eqnarray}}
       \newcommand{\ea}{\end{eqnarray}}
\newcommand{\ban}{\begin{eqnarray*}}
\newcommand{\ean}{\end{eqnarray*}}
\newcommand{\Bbb}[1]{\mbox{\BBb #1}}

\newcommand{\NN}{\Bbb{N}}
\newcommand{\ZZ}{\Bbb{Z}}

\newcommand{\RR} {\Bbb{R}}
\newcommand{\TT} {\Bbb{T}}

\newcommand{\EE} {\Bbb{E}}

 \newcommand{\Pf}{\noindent {\bf Proof:} }

\newcommand{\sect}[1]{\section{#1} \setcounter{equation}{0}}
\newtheorem{theorem}{Theorem}[section]
\newtheorem{theo}{Theorem}[section]
\newtheorem{remark}[theo]{Remark}
\newtheorem{rmk}[theo]{Remark}
\newtheorem{defn}{Definition}[section]
\newtheorem{example}{Example}[section]

\newtheorem{lemma}[theo]{Lemma}
\newtheorem{prop}[theo]{Proposition}
\newtheorem{coro}[theo]{Corollary}

\newtheorem{problem}[theo]{Problem}
\begin{document}
\title{Convergence and the Length Spectrum}

\author{C. Sormani
\thanks {Partially supported by NSF Grant \# DMS-0102279
and a PSC CUNY Award}}

\date{}
\maketitle

\noindent{\bf Abstract:} 
The author defines and analyzes 
the $1/k$ length spectra, $L_{1/k}(M)$, whose union, 
over all $k\in \NN$ is the classical length spectrum.  
These new length spectra are shown to
converge in the sense that 
$\lim_{i\to\infty} L_{1/k}(M_i) \subset \{0\}\cup L_{1/k}(M)$
as $M_i\to M$ in the Gromov-Hausdorff sense.   
Energy methods are introduced to estimate the shortest 
element of $L_{1/k}$, 
as well as a concept called the minimizing index which may be used
to estimate the length of the shortest closed geodesic of a simply
connected manifold in any dimension.
A number of gap theorems
are proven, including one for manifolds, $M^n$, with $Ricci\ge (n-1)$
and volume close to $Vol(S^n)$.   
Many results in this paper hold on compact length spaces in addition to
Riemannian manifolds.

\sect{Introduction} \label{intro} 

Recall that a compact length space is a metric space such that
every pair of points is joined by a length minimizing rectifiable 
curve whose length is the distance between the two points.  
The simplest example of such a space is a Riemannian manifold.
A ``geodesic'' in such a space is a locally length minimizing curve.

A closed geodesic is a map $\gamma: S^1 \to M$ which is locally minimizing
around every point in $S^1$.  This extends the concept of a smoothly 
closed geodesic in a manifold.  (c.f. \cite{BBI}{SoWei})  We shall
assume throughout that all geodesics are parametrized proportional to 
arclength with speed $L/(2\pi)$.  The length spectrum, $L(M)$,
of a length space, $M$, is
the set of lengths of closed geodesics.  These definitions are
just extensions of the clasical definitions on Riemannian manifolds.

The length spectrum is not continuous with respect to deformations
of the manifold.  When a sequence of spaces, $M_i$, converges 
in the Gromov-Hausdorff sense [Defn~\ref{defnGH}]
to a space, $M$, it may have closed geodesics $\gamma_i$ converging to
a curve which is not a closed geodesic.  That is, there could be a 
``disappearing length'':
$\exists L_0=\lim_{i\to\infty} L_i \in L(M_i)$ such that $L_0\notin L(M)$.

%\psdraft
\begin{figure}[htbp]
\includegraphics[width=6in ]{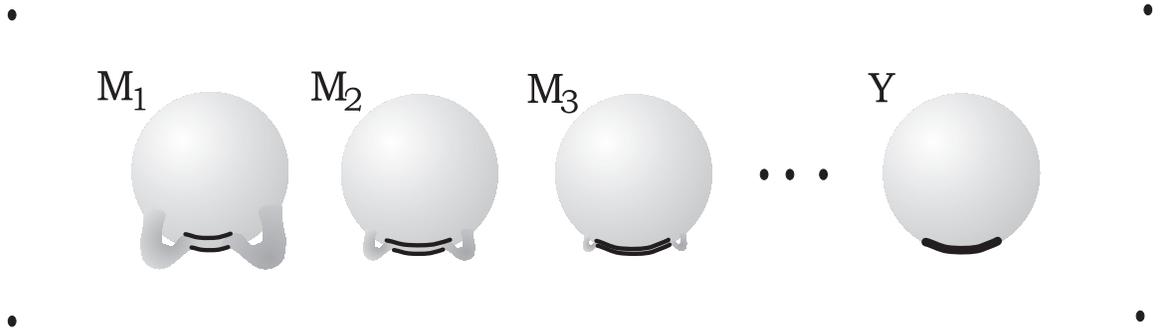}
\caption{$M_i \to Y$ in the Gromov Hausdorff sense but a length disappears. }
\label{tosegment}
\end{figure}
%\psfull

In particular, we have this situation in 
Figure~\ref{tosegment}. Here the sequence of surfaces, $M_i$, with  
increasingly small pairs of handles
converges in the Gromov-Hausdorff sense to a standard sphere, $Y$.
Notice how the closed geodesics which pass through both handles
converge to a geodesic segment but not to a closed geodesic.
The lengths, $L_i$, of these closed geodesics converge
to $L_0=\pi/3$ which is not in the length spectrum of the
sphere.  In fact the shortest closed geodesic in $S^2$ has
length $2\pi$.
For details see Example~\ref{handles}.

In Example~\ref{torusH}, we will
examine the length spectrum of a flat torus created by taking the
isometric product of a circle with a small circle.  As the smaller
circle's diameter approaches $0$, we say the sequence of tori ``collapses''
in the Gromov-Hausdorff sense to a circle, $S^1$.  
The length spectrum of
the limit space, $S^1$, is just $\{n\pi: n\in \NN\}$, yet the length spectra
of the collapsing tori is becoming an increasingly dense collection
of points in $[0,\infty)$.  Thus we have quite a large collection of
disappearing lengths!

Both of these examples will be described in full detail in the first
section [Example~\ref{torusH}] and [Example~\ref{handles}], after
we have given the rigourous definition of Gromov-Hausdorff convergence.

\begin{figure}[htbp]
\includegraphics[width=6in ]{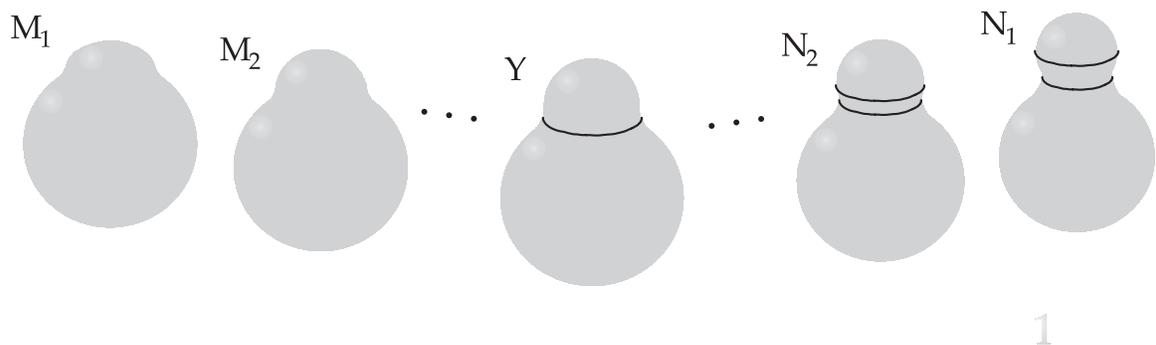}
\caption{ The geodesic in $Y$ is suddenly appearing as a limit of 
the $M_i$ but not as a limit of the $N_i$. }
\label{boobah}
\end{figure}

It is also possible that there is a ``suddenly appearing length'':
$L_0 \in L(M)$ such that no sequence $L_i \in  L(M_i)$
converges to $L$.  This occurs even when $M_i$
converges to $M$ in the $C^4$ sense as can be seen in Figure~\ref{boobah}.

%\psdraft
%\psfull

In this paper, we define a new collection length spectra,
$L_{1/k}(M)$,  [Defn~\ref{defL1/k}] such that
\be \label{bigcupintro}
\bigcup_{k\in \NN} L_{1/k}(M)= L(M) \textrm{ [Theorem~\ref{bigcup}]}.
\ee
While any collection of length spectra satisfying (\ref{bigcupintro})
would have to incoporate the sudden appearances observed in
Figure~\ref{boobah}, we do prove in Theorem~\ref{nodis} that
\be \label{nodisintro}
\lim_{k\to \infty} L_{1/k}(M_i) \subset \{0\} \cup L(M),
\ee
when $M_i$ converge to $M$ in the Gromov-Hausdorff sense.

Throughout the paper we survey past results and techniques used
to study the length spectrum, we relate them to the new $1/k$ length
spectra and we suggest new directions of research.  As many of the
proposed problems in this paper are at a level a graduate student
should be able to handle, we have presented this paper in a manner
that should easily be read by a student.

In Section~\ref{background} we give the necessary background on
Gromov-Hausdorff convergence, completely describing how the 
spheres with tiny handles and the collapsing tori converge
[Examples~\ref{handles} and~\ref{torusH}].  We also present
ellipsoids which converge to a singular doubled disk Example~\ref{DD}.
Readers who are just interested in the $1/k$ length spectra and not
their convergence properties may skip Section~\ref{background}
and easily read everything except Sections~\ref{nodis} and~\ref{almrig}.

In Section~\ref{sectL1/k}, we introduce $1/k$ geodesics, which
are geodesics that minimize on any subsegment whose length is
$1/k$ of the total length [Defn~\ref{1overk}].  The set
of lengths of such geodesics is $L_{1/k}$ [Defn~\ref{defL1/k}]
and we prove Theorem~\ref{bigcup}.  We also relate $L_{1/k}(M)$
to the diameter and injectivity radius of the space [Lemma~\ref{diam}
and Lemma~\ref{injLk}].  Using these results we descibe the $1/k$
length spectra of a sphere and collapsing tori [Examples~\ref{sphere0}
and~\ref{torus0}].  

In Section~\ref{sectL1/k} we also complete a study of closed geodesics.
We define the minimizing index of a geodesic as the smallest $k$ which can be
used to classify it as a $1/k$ geodesic [Defn~\ref{minind}].
Then we define the injectivity radius of a geodesic [Defn~\ref{injgamma}]
and relate it to the minimizing index  [Lemma~\ref{ellipsoidlem}].
We discuss iterated geodesics [Lemma~\ref{w2}] and a particularly 
illustrative example of the equator of an ellipsoid close to a doubled disk 
[Example~\ref{ellipsoid0}].

In Section~\ref{sectcov}
we prove that the covering spectrum defined in \cite{SoWei}
is a subset of $L_{1/2}$ [Theorem~\ref{covspec}].
Recall that in \cite{SoWei}, the $CovSpec(M)\cup \{0\}$
was proven to be continuous with respect to Gromov-Hausdorff
convergence of the manifold.  In other words, there is no
sudden appearance of elements in the Covering Spectrum
as described in Figure~\ref{boobah}.  We could not expect to
prove such a strong theorem for $L_{1/k}(M)$ because we include
all lengths of $L(M)$ in one of the them [Theorem~\ref{bigcup}].
In fact the suddenly appearing geodesic in Figure~\ref{boobah}
is a $1/2$ geodesic and an element of $L_{1/2}$.  This justifies the
lack of an equality in Theorem~\ref{nodis}.

In Section~\ref{secthalf} we turn to a study of the systole
of a manifold.  Since the systole is the length of the shortest
noncontractible curve it is an element of $L_{1/2}$ [Lemma~\ref{w1}].
We survey
past estimates relating the systole to the volume and diameter of
a manifold and extend them to estimates on $\min L_{1/2}$.
It should be noted that some of these estimates are only achieved
on singular manifolds, so the extension of all concepts to 
compact length spaces in this paper is further justified.

In Section~\ref{sectmin} we estimate the length of the shortest
closed geodesic in a compact length space, $\min L(M)$.
First we define the minimizing index, $minind(M)$, of a space 
and then prove that
$\min L(M)\le minind(M)diam(M)$ [Defn~\ref{minindM} and 
Theorem~\ref{minindMd}].
We also provide a lower bound on $\min L(M)$ [Theorem~\ref{minindMinj}].
We close with an application of an old result of Klingenberg
to show that the minimizing index of a manifold without conjugate points
is $1/2$ [Corollary~\ref{klingcor}].

In Section~\ref{sectnodis}, we finally prove the convergence
theorem mentioned above.  We conclude that if there is a sequence
of spaces with a disappearing length in the limit, as in 
Figure~\ref{tosegment}, then the geodesics that disappear must
have minimizing index diverging to infinity [Corollary~\ref{minindtoinfty}].
Theorem~\ref{nodis} also immediately implies that $L_{1/k}(M) \cup \{0\}$
is compact [Corollary~\ref{L1/kcompact}].   We discuss the
convergence of $L_{1/k}(M_i)$ for the collapsing tori, the flattening
ellipsoids and  a new example converging to a hexagonal region
[Examples~\ref{torus1}, ~\ref{ellipsoid1} and~\ref{hexagon1}].
We remark on Bangert's Theorem [Remark~\ref{bangertrmk}].

In Section~\ref{almrig} we rephrase Theorem~\ref{nodis} as
a gap thoerem [Theorem~\ref{gapthmk}] and then prove a number of 
gap theorems.  For example, we use Colding's sphere stability theorem
to prove:

\begin{theorem}\label{volgap} 
There exists a function $\Psi:\RR^+\times \NN \times \NN \to \RR^+$ such
that $\lim_{\delta \to 0}\Psi(\delta, k,n)=0$ such that
if $N^n$ is a compact Riemannian manifold with
\be
Vol(N^n) \ge Vol(S^n)-\delta
\ee
and $Ricci(N^n) \ge (n-1)$ then
\be
L_{1/(2k)}(M^n) \subset [0,\epsilon_k)\cup (2\pi-\epsilon_k, 2\pi+\epsilon_k) 
\cup \cdots (2k\pi-\epsilon_k, 2k\pi+\epsilon_k)
\ee
for $\epsilon_k=\Psi(\delta,k,n)$.
\end{theorem}

In light of Example~\ref{football}, we propose that the length spectra
on these manifolds do not converge [Problem~\ref{volgapdis}.
We obtain similar results for Riemannian manifolds
with $Ricci\ge (n-1)$ and $rad(M^n)$ close to $rad(S^n)$ 
and for $M^n$ with first
Betti number equal to $n-1$ and Ricci curvature that is almost
nonnegative [Theorems~\ref{radgap} and~\ref{torusgap}].
Similar gap theorems also exist for $M^n$ which are
almost isotopic and have a uniform lower bound on their Ricci
curvature [Remark~\ref{cosmos}].

In Section~\ref{sectopen}, we introduce openly $1/k$ geodesics
which are shown to be
uniquely defined on Riemannian manifolds by any collection of
evenly spaced points [Defn~\ref{defopen} and Lemma~\ref{oldunique}].
We then define the openly $1/k$ length spectra and extend most
of the results and definitions of the previous sections to this setting.  
There are no openly $1/2$ geodesics [Lemma~\ref{back}], so the
results in Sections~\ref{sectcov} and~\ref{secthalf}  do not apply.
Theorem~\ref{nodis} doesn't extend well either due to the open
nature of the Defn~\ref{defopen} [Theorem~\ref{opnodis} and 
Example~\ref{ellipsoidopen}].  Otherwise the results carry over.
We close with a discussion of the borderline case of a $1/k$ geodesic
which is not an openly $1/k$ geodesic on manifolds.

Section~\ref{sectenergy} extends the theory of geodesics on  as critical
points of the energy function on the loop space to openly $1/k$
geodesics.  
First we review the piecewise geodesic version of the
theory and then prove Theorem~\ref{openenergy} which identifies
openly $1/k$ geodesics on a convex compact Riemannian manifold with boundary
 with ``rotating'' smooth critical points of a
uniform energy function on $k$-fold product, $(M)^k$, of the space.
We explicitly demonstrate a few examples and then discuss
why this thoery does not extend well to $1/k$ geodesics and
nonsmooth spaces.  Nevertheless it can be used to estimate the minimizing
index of a Riemannian manifold and to determine whether a
convex Riemannian manifold with boundary has any closed geodesics.

Section~\ref{sectproblems} concludes the paper with a list of open
problems most of which should be on the level of a graduate student.

Additional gap theorems related to sectional curvature will appear
in future papers along with a survey of sectional curvature results on the
length spectrum .

The author is grateful to Guofang Wei (UCSB) for many fruitful suggestions,
to Wolfgang Ziller (U Penn) for guidance through the literature,
to Steve Zelditch (JHU) for requesting a thorough extension
of the convergence results in \cite{SoWei}, and
to Carolyn Gordon (Dartmouth) for her excellent advise in all things.

\sect{Background} \label{background}

Here we provide the necessary background on Gromov-Hausdorff
convergence.  Readers
who are only interested in studying the $1/k$ length spectrum on
a fixed Riemannian manifold may skip to Section~\ref{sectL1/k}.
Essentially all the material here appeared in \cite{G3} and can
also be studied in \cite{BBI}.

For those readers who have studied $C^k$ convergence of
manifolds, keep in mind that Gromov has proven that
if a sequence of compact manifolds $M_i$
converges to $M$ in the $C^k$ sense then they also converge in the 
Gromov-Hausdorff sense.  While $C^k$ convergence requires that the manifolds
be diffeomorphic, $GH$ convergence doesn't even require that they have
the same dimension.  In fact the spaces need only be compact metric spaces.

We begin with an older concept, the Hausdorff distance between sets.

\begin{defn}\label{defnH}
Given two compact subsets $A,B$ in a metric space $Z$, we can define
the Hausdorff distance as follows:
\be
d_H^Z(A,B) = \inf\{r: \,\, A \subset T_r(B) \textrm{ and } B\subset T_r(A) \},
\ee
where $T_r(X)$ is the tubular neighborhood around $X$ of radius $r$:
\be
T_r(X)=\{z\in Z: \,\, \exists x_z\in X \,\,s.t.\,\, d(x_z,z)<r\}.
\ee
\end{defn}

The surfaces $M_i$ in Figure~\ref{tosegment} would converge
in the Hausdorff sense to the standard sphere, $Y$, as subsets of 
$Z=\EE^3$, if they were superimposed.
One need only take the radius of the tubular neighborhood large
enough to capture the tiny handles.  In this respect the Hausdorff
distance is blind to the topology of the sets it compares.  

Hausdorff distance
is also blind to the dimensions of the sets: it can easily be
seen that $A_0=[0,1]\times \{0\}\subset \EE^2$ and 
$A_r=[0,1]\times [-r,r] \subset \EE^2$ satisfy
$d_H(A_0,A_r)\le r$.  On the other hand, for small $r$,
one can see how it makes sense that $A_r$ could be
thought of as close to $A_0$.  To quote Cheeger, they
look very similar ``to the naked eye''.

Gromov extended this concept to compact metric spaces, providing
us with a metric between spaces that is also blind to dimension
and topology, but captures the idea that the spaces are close
in some blurry sense.  \cite{G3}

Before we define the Gromov-Hausdorff distance between metric spaces,
recall that $f:X\to Z$ is an isometric embedding if it is one-to-one and
$d_X(x_1,x_2)=d_Z(f(x_1), f(x_2))$ for all $x_1, x_2 \in X$.  The sphere
sitting inside Euclidean space is not isometrically embedded because 
the distances on the sphere are measure intrinsically (the poles are
a distance $\pi$ apart in $S^2$).  A plane is isometrically embedded
in Euclidean space because it is totally geodesic.

\begin{defn} [Gromov]\label{defnGH} 
The Gromov-Hausdorff distance between two compact metric spaces
$X$ and $Y$ is defined as follows:
\be
d_{GH}(X,Y)=\inf \{ d_H^Z(f(X),g(Y)): Z, f:X\to Z, g:Y\to Z \}
\ee
where the set runs through any metric space, $Z$,
and any isometric embeddings $f:X\to Z$ and $g:Y \to Z$.
\end{defn}

It is an easy exercise to prove the following two lemmas:

\begin{lemma} [Gromov]
\be
d_{GH}(X,Y) \le diam(X) + diam(Y)
\ee
\end{lemma}

\begin{lemma} [Gromov] \label{GHdiam}
If $X_i$ converge to $X$ in the Gromov-Hausdorff sense,
$d_{GH}(X_i,Y) \to 0$ then $diam(X_i) \to diam(X)$.
\end{lemma}

Gromov proved that both
the space of compact metric spaces and  the space of
compact length spaces are complete with respect to $d_{GH}$.  
In particular, he proved the difficult theorem that
the Gromov-Hausdorff limit of a compact length space is
a compact length space.  \cite{G3}

We can now give the details of the sequence of tori collapsing 
to a circle with disappearing lengths that was
mentioned in the introduction.

\begin{example}\label{torusH} {\em
Let $M_j=S^1_{\pi}\times S^1_{\pi/j}$ be a flat torus
formed by taking the isometric product of a circle of diamter
$\pi$ with a circle of diameter $1/j$.  
Note that as $j$ diverges to infinity,
$M_j$ converges in the Gromov-Hausdorff sense to $S^1_\pi$.
This can be seen by taking $Z=M_j$ itself and isometrically
embedding $S^1_\pi$ as $S^1_\pi \times \{0\}\subset M_j$, so
\be
d_{GH}(M_j, S^1_\pi)\le d_H^{M_j}(M_j, S^1_\pi \times \{0\}) \le \pi/j.
\ee

It is well known that the length spectrum of $M_j$ is the collection
of distances between lattice points $(2\pi a, 2\pi b/j)$ where
$a,b \in \ZZ$.  Thus
\be
L(M_j)=\{\sqrt{(2\pi a)^2 + (2\pi b/j)^2}:\, a,b \in \NN \} 
\, \cup \,  \{\pi, \pi/k\}.
\ee
Notice how this length spectrum becomes increasingly dense as
$j$ goes to infinity, so that for any $N$ we get
\be
L(M_j)\cap [0, N] \to [0,N] \textrm{ in the Hausdorff sense.}
\ee
In particular, there are lengths $l_j \in L(M_j)$ such that
$l_j \to \pi$ even though $\pi$ is not in the length
spectrum of $S^1$.  
}
\end{example}

The sequence of surfaces in Figure~\ref{tosegment} are trickier to
deal with as they are not easily isometrically embedded into a 
common space and even the choice of a metric space $Z$ for each
$M_j$ is not obvious.  

We first recall Gromov's concept of an $r$ net on a metric space.
A set $N \subset X$ is an $r$ net if $X\subset T_r(N)$.  
It is clear that $d_H(N,X)\le r$.  When $N$ is a finite
collection of points then it is a finite net.
Let $N_X(r)$ be the minimum cardinality of all
$r$ nets in $X$.

Gromov's famous compactness theorem states that a class of 
compact metric spaces, $\{X\}$, with a uniform bound on
$N_X(r)\le N(r)$
is precompact with respect to the Gromov-Hausdorff metric.
In particular, the class of
complete manifolds with $Ricci \ge -K$, $dim=n$ and $diam\le D$ is
precompact \cite{G3}\cite{Bi}.  The limits of the sequences are
compact length spaces. \cite{G3}

If one considers an r net, $N \subset X$, and endows it with the
restricted metric from $X$, then it may not be a length space.
However, it is a metric space such that
$d_{GH}(X, N) \le d_H^X(X,N)\le r$.  Using the triangle inequality,
one than sees that $d_{GH}(X,Y)\le 2r$ if both spaces have
isometric $r$ nets.  We now use this technique to prove the convergence
of the surfaces in Figure~\ref{tosegment}.

\begin{example} \label{handles}
The surfaces in 
Figure~\ref{tosegment} converge to the standard sphere
in the Gromov-Hausdorff sense.  Let us suppose that $M_j$
with it's handles removed is isometric to a
a standard sphere with two disks of radius $1/j$ removed
and that the diamater of the handles is $<4/j$.  
Now lets form a finite $100/j$ net on $M_j$ such that for any
pair of points in the net, the minimizing geodesic between
them does not hit either handle.  Since the points in the net
aren't on the handles, they correspond isometrically
to specific points on $Y=S^2$.  That is we have a metric space
$N_j$, the net, such that $N_j$ isometrically embeds into
both $Y$ and $M_j$ and such that 
\be
d_H^{M_j}(N_j, M_j)\le 100/j \textrm{ and }
d_H^{Y}(N_j, Y)\le 100/j.
\ee
Thus
\be
d_{GH}(M_j, Y)\le d_{GH}(M_j,N_j)+d_{GH}(N_j,Y)\le 200/j
\ee
and we see that $M_j$ converges to $Y$ is the Gromov-Hausdorff
sense.
\end{example}

Note that it is not necessary to find an isometric pair of $r$
nets in two compact metric spaces $X_1$ and $X_2$ to prove they
are close in the Gromov-Hausdorff sense.  It suffices to find
a pair of ``almost isometric'' nets $N_1$ and $N_2$.  
Gromov has proven that
if both nets have the same cardinality and one can set up a
bijection between them: $f:N_1 \to N_2$ such that
\be \label{CkGH}
sup \{ |d_{N_1}(x_1, x_2)-d_{N_2}(f(x_1), f(x_2))|: x_1, x_2\in N_1\} <\epsilon
\ee
then one can show $d_{GH}(N_1, N_2) < 2\epsilon$ [Gr], cf [BBI. Cor 7.3.28].
So  in that case $d_{GH}(X_1, X_2)<2r+\epsilon$.

Using (\ref{CkGH}) we see that when $X_i \to X$ in the $C^k$ sense
then they also converge in the Gromov-Hausdorff sense.  
For example, Figure~\ref{boobah} has a smoothly converging
sequence of compact Riemannian manifolds.  They are all diffeomoerpic
to the sphere with metrics $g$ that converge smoothly, $C^k$, to
the limit space.  The diffeomorphisms are almost isomorphisms
in the sense described in (\ref{CkGH}) without even requiring
the use of finite nets.

In fact one need only find $f_i:X_i \to X$ which are
$\epsilon_i$ almost distance preserving,
\be
|d_X(f_i(a), f_i(b))-d_{X_i}(a,b)|<\epsilon_i
\ee
 and $\epsilon_i$
almost onto, $T_\epsilon(f_i(X_i))\supset X$, to prove
that $X_i$ converge to $X$ in the Gromov Hausdorff sense.

\begin{example} \label{DD}
Let $M_{c}$ be an ellipsoid 
\be
(x)^2+(y)^2 +(z/c)^2 =1.
\ee
If we take $c_j\to\infty$, then $M_j=M_{c_j}$
converges to the doubled disk, $Y$, in the 
Gromov Hausdorff sense.  

More precisely, $Y$, is two flat disks of radius 1 glued
together along the circle, so that the distance between points
on a common disk is the usual Euclidean distance and the
distance between points $x$ and $y$ on
different disks is:
\be
d_{M_\infty}(x, y)=\inf_{z\in S^1}(|x-z|+|y-z|).
\ee
Gluing is significantly more complicated when the shapes aren't
convex (c.f.\cite{BBI}).  

To prove that $M_j$ converge to $Y$ in the Gromov-Hausdorff 
sense, we just employ the maps $f_i:M_i \to M$, defined to
be $f_i(x,y,z)=(x,y, sgn(z))$, where $sgn(z)$ is just used to
indicate whether we are on the upper or lower disk.  Naturally
the edge where $z=0$ doesn't need a sign.
\end{example}

Example~\ref{torusH} is said to be ``collapsing'' because the
dimension of the limit is less than the dimension of the
sequence.  
Example~\ref{handles} is considered to be ``noncollapsing'' 
because the dimension of the manifolds in the sequence is
the same as the dimension of the limit space.  On the other hand,
the injectivity radius is decreasing to $0$ in this example.

The following definition is a simple extension of a Riemannian
injectivity radius.

\begin{defn} \label{injdef}
The injectivity radius of a compact length space, $M$,
is
\be
injrad_x=\sup \{ t: \textrm{ any geodesic segment
of length } t \textrm{ is minimal.}\}
\ee
\end{defn}

\begin{example} \label{Hawaii}
The Hawaiian Earring, a compact length space consisting of a collection
of circles of radius $1/j$ for each $j\in \NN$ all joined at 
a common  point,
has an injectivity radius equal to 0.
It is also known to have no universal cover (c.f. \cite{Sp}).
That is, there is no covering space which covers all the other
covering spaces.  
%It's length spectrum is $\QQ \cup (0, \infty)$.
%wrong it is $\RR^+$.
\end{example}

For completeness of exposition,
we now review the fact that there are no disappearing lengths
when the sequence of compact length spaces 
has a uniform lower bound on injectivity
radius.

\begin{prop} \label{propinj}
Suppose $M_j$ are compact length spaces with a common
positive lower bound, $i_0$, on their injectivity radiui,
and $M_j$ converge to $Y$ in the Gromov-Hausdorff metric,
then for any $R>0$ we have the following Hausdorff limit:
\be
L(M_j)\cap [0,R] \to L_\infty \subset L(Y)\cap [0,R].
\ee
That is, for all $\epsilon>0$, $\exists N_{\epsilon,R}\in \NN$
such that
\be
L(M_j)\cap [0,R] \subset T_\epsilon(L(Y)).
\ee
\end{prop}

Note that when manifolds converge smoothly, they do have
a common lower bound on their injectivity radius.  So
Figure~\ref{boobah} demonstrates that one still might
have suddenly appearing geodesics in this case.  

One can
see that $N_{\epsilon,R}$ depends on $R$, just by examining
a sequence of circles of radius $r_j\to \pi$
converging to the standard circle.  The errors
accumulate as you wrap repeatedly around the same geodesic.
So one needs a common upper bound, $R$, on the length of
the geodesic to get a common rate of convergence.

\Pf
 Let $\gamma_i$ be the geodesics of length $L_i \to L_\infty$.
For $i$ sufficiently large, $L_i \in [2i_0, 2L_\infty]$.  
Thus the $\gamma_i:S^1 \to M_i$ are equicontinuous.
By the Grove-Petersen Arzela-Ascoli Theorem, a subsequence
of the $\gamma_i$ converge to $\gamma_\infty:S^1\to Y$,
of length $L_\infty$. \cite{GrPet}
\qed

One may consider Proposition~\ref{propinj} as a kind of
semicontinuity of the length spectrum with respect
to Gromov-Hausdorff convergence.  One of the goals of
this paper is to prove a similar convergence theorem
for spaces without a common lower bound on injectivity
radius.

%=======================================================

\sect{1/k Geodesics}  \label{sectL1/k}

We now introduce a new length spectrum,
$L_{1/k}$, which we will later prove has a strong
relationship with Gromov-Hausdorff convergence.
Here we focus on the properties of this new
concept on a fixed compact length space or Riemannian 
manifold and it's relationship with the traditional
length spectrum.

\begin{defn} \label{1overk} 
A $1/k$ geodesic is a closed geodesic, $\gamma:S^1 \to M$,
which is minimizing on all subsegments of length $L/k$
where $L=Length(\gamma)$:
\be \label{1overk1}
d(\gamma(t), \gamma(t+2\pi/k))=L_\gamma([t, t+2\pi/k])=
L/k \qquad \forall t\in S^1.
\ee
\end{defn}

\begin{defn} \label{defL1/k}
Let the $1/k$ Length Spectrum, $L_{1/k}(M)\subset L(M)$,
be the set of lengths of $1/k$ geodesics.
\end{defn}

The following lemma is immediate:

\begin{lemma} \label{w0}
$L_{1/k}(M)\subset L_{1/(k+1)}(M)$.
\end{lemma}

\begin{defn} \label{minindex}\label{minind}
The minimizing index, $minind(\gamma)$,  
of a geodesic, $\gamma$,
is the smallest $k \in \NN$ such
that the geodesic is a $1/k$ geodesic.
\end{defn}

\begin{theorem} \label{bigcup}
Any closed geodesic is a $1/k$ geodesic
for a sufficiently large number $k$.  So
\be
\bigcup_{k=1}^{\infty} L_{1/k}(M)= L(M).
\ee 
\end{theorem}

\Pf
If $\gamma:S^1\to M$ is a closed geodesic
then for all $t$, there exists $\epsilon_t>0$
such that $\gamma$ is minimizing from
$t-\epsilon_t$ to $t+\epsilon_t$.  The intervals
$(t-\epsilon_t, t+\epsilon_t)$ form an open cover
of $S^1$ and since $S^1$ is compact, there is
a finite subcover and a Lebesgue number, $\rho>0$, for
this cover.  The $\gamma$ is a $1/k$ geodesic
for any $k>2\pi/\rho$.
\qed

\begin{lemma} \label{diam}  
If $diam(M) \le D$ then $minind(\gamma)\ge L(\gamma)/D$ and 
$L_{1/k}(M) \subset (0, Dk]$.
\end{lemma}

\Pf
If $minind(\gamma)=k$, then
$\gamma$, must be minimizing on 
segments of length
$L(\gamma)/k$.  So $L(\gamma)/k \le D$.
\qed

Recall Definition~\ref{injdef} of the injectivity radius..
It is easy to see that $L(M) \subset [2injrad(M), \infty)$.
The injectivity radius also provides the following
useful relationship between $L(M)$ and
$L_{1/k}(M)$.

\begin{lemma} \label{injLk} 
If $M$ has $injrad(M)\ge i_0>0$ then for all $L_0>0$,
\be
L(M)\cap [0,L_0]= L_{1/k}(M) \cap [0, L_0]  
\textrm{ where } k\ge L_0/i_0.
\ee
\end{lemma} 

\Pf
Let $\gamma:S^1 \to M$ be any closed geodesic with $L(\gamma)\le L_0$.  
Then any segment
from $\gamma(t)$ to $\gamma(t+L/(2\pi k))$ has length $\le L_0/k \le i_0$.
Thus it is minimizing on this interval and $\gamma$ is a $1/k$ geodesic.
\qed 

This estimate is only sharp when the injectivity radius of the manifold
is achieved by a pair of points on the geodesic.  There is no reason 
that a distant pair of cut points or a cut point perpendicular to the
geodesic should affect its minimizing index.  In Example~\ref{torus0}
below we will see that in a thin torus, $S^1_\delta\times S^1_\pi$,
there is a $1/2$ geodesic of length $2\pi$ no matter how small
the injectivity radius, $\delta$, of the torus is.  For this reason we
make the following new definition.

\begin{defn} \label{injgamma}
The injectivity radius, $injrad(\gamma)$ of a closed 
geodesic $\gamma$, is
\be
injrad(\gamma)=\inf\{h_t: t\in S^1\}
\ee
where
\be
h_t= \sup \{h: \gamma \textrm{ is minimizing on [t, t+h].} \}  
\ee
\end{defn}

The following lemma is easily follows from Definitions~\ref{1overk} 
and~\ref{injgamma}.

\begin{lemma} \label{injgammalem} \label{ellipsoidlem}
A closed geodesic $\gamma$ of length $L$ satisfies
\be
\frac{L}{\, minind(\gamma) \,} \, \le \, injrad(\gamma) \, < \,
\frac{L}{\, (minind(\gamma)-1) \, } .
\ee
\end{lemma}

Recall that a prime geodesic is a closed geodesic
whose period is $2\pi$.  
All closed geodesics
are either prime geodesics or iterated geodesics
of the form $\gamma_n(t)=\gamma_1(nt)$ where $\gamma_1$
is a prime geodesic, and $n\in \NN$.  

\begin{lemma}  \label{w2}  
If $\gamma_1: S^1 \to M$ is a $1/k$ closed geodesic and 
$\gamma_n:S^1\to M$ is defined by $\gamma_n(t)=\gamma_1(nt)$, 
then $\gamma_n$ is an
$1/(kn)$ geodesic.  In fact
%w2 
\be 
minind(\gamma_n) \in [n(minind(\gamma_1)-1), n\,minind(\gamma_1)) ]
\cap [2n, \infty).
\ee
\end{lemma}

\Pf
Let $L$ be the length of $\gamma_1$ and $k=minind(\gamma_1)$.
Then $nL=L(\gamma_n)$ and
\be
d(\gamma_n(t), \gamma_n(t+(2\pi/(kn)))=d(\gamma_1(nt), \gamma_1(nt+2\pi/(k)))
=L/k=L(\gamma_n)/(nk),
\ee
which implies that  $minind(\gamma_n)\le (nk)$.
On the other hand, suppose $j=minind(\gamma_n)$, then
\be
nL/j \, =\, d(\gamma_n(t), \gamma_n(t+(2\pi/j))=
d(\gamma_1(nt), \gamma_1(nt+2\pi n/j))
\ee
So $injrad(\gamma_1) \ge n L(\gamma_1)/j$ and applying 
Lemma~\ref{injgammalem} we get
\be
L(\gamma_1)/(minind(\gamma_1)-1) \ge injrad(\gamma_1) \ge n L(\gamma_1)/j
\ee 
which implies $minind(\gamma_1)-1 \le (j/n) \le minind(\gamma_n)/n $.

Our final consideration is that $injrad(\gamma_1)\le L(\gamma_1)/2$,
so $minind(\gamma_n)\ge 2n$.
\qed

We can now apply these lemmas to examine some examples.

\begin{example}\label{sphere0}
Suppose $S^2$ is the standard sphere.  It
is well known that all its prime
closed geodesics have length $2\pi$.
These geodesics can easily be seen to
be $1/2$ geodesics.  By Lemma~\ref{w2},
we then have
\be
2k\pi \in L_{1/2k}(S^2)
\ee
and by Lemma~\ref{w0},
\be
\{2\pi, 4\pi, ... 2k\pi\} \subset L_{1/2k}(S^2)\subset L_{1/(2k+1)}(S^2).
\ee
This also follows directly from Lemma~\ref{injLk}.
On the other hand, by Lemma~\ref{diam},
\be
L_{1/j}(S^2)\subset L(S^2)\cap (0, j\pi]
\ee
which gives  
\be
L_{1/2k}(S^2)= L_{1/(2k+1)}(S^2)=\{2\pi, 4\pi, ... 2k\pi\}.
\ee
\end{example}

\begin{example}\label{torus0}  {\em
Let $M_j=S^1_{\pi}\times S^1_{\pi/j}$ be a flat torus
from Example~\ref{torusH}.  The closed geodesics of the torus
are of the form
\be
\gamma(t)=\left( (a t +x_0) mod 2\pi , (b t/j +y_0) mod 2\pi/j \right), 
\ee
where $a, b \in \ZZ$ and $x_0, ky_0 \in [0,2\pi]$.
It is minimizing until $|a|t=\pi$ or $|b|t/j=\pi/j$, that
is until $t = \min \{\pi/|a|, \pi/|b|\}$.  So it's minimizing index
is 
\be
minind(\gamma)=\max \{2|a|, 2|b|\}.
\ee

Note that $\gamma$ is a prime geodesic whenever $a$ and $b$ are
relatively prime or $ab=0$ and one of them has absolute value 1.
In particular the geodesic with $b=a+1$ is a prime geodesic.
{\em So for any natural number, k, there is a prime geodesic
in the torus of with minimizing index $=2k$.}

The length of our arbitrary geodesic, $\gamma$, is
\be
L(\gamma)=\sqrt{(2\pi a)^2 + (2\pi b/k)^2 \, }
\ee
So, skipping the trivial geodesic, we have
\be
L_{1/(2k)}(M_j)=\{\sqrt{(2\pi a)^2 
+ (2\pi b/j)^2 \,}: a,b = 0,1,2,...k \}\setminus \{0\}.
\ee
and $L_{1/(2k-1)}(M_j)=L_{1/(2k)}(M_j)$. }
\end{example}

Recall Lemma~\ref{w2} which states that the nth iterate of $1/k$
geodesic is a $1/(nk)$ geodesic.  This does not mean its minimizing
index is $nk$.  In fact in Example~\ref{ellipsoid0} we will
demonstrate that the minimizing index of an iterated geodesic
may take on any natural number allowed in the Lemma.  

\begin{example} \label{ellipsoid0}
We claim that for any $k\in \NN$ there exists $c_k\in (0,1]$ such that the
the ellipsoid, $M(c_k)$: 
\be
(x)^2+(y)^2 +(z/c_k)^2 =1.
\ee
has a prime geodesic 
$\gamma_{c_k}(t)=(cos(t), sin(t),0)$ whose minimizing index is $k+1$.
Furthermore, for any $n\in \NN$, there exists $c_{n,k} \in (0,1]$
such that $\gamma_c(nt)$ has minimizing index equal to $k+n$.

The brute force proof of the claim is to use the recent work of 
Itoh and Kiyohara to explicitly determine the cut locus
of the points on this geodesic \cite{ItKi}.  One sees 
that $injrad(\gamma_c)$ varies continuously with $c$ taking on 
all values in $(0,\pi)$. 
The claim then follows by applying Lemma~\ref{ellipsoidlem}
%Itoh, Jin-ichi(J-KUMAED); Kiyohara, Kazuyoshi(J-OKAYS) 
%The cut loci and the conjugate loci on ellipsoids. (English. English summary) 
%Manuscripta Math. 114 (2004), no. 2, 247--264.
\end{example}

Lemma~\ref{diam} implies that the minimizing
index of a closed geodesic, $\gamma: S^1\to M$ satisfies
\be
minind(\gamma) \ge L(\gamma)/Diam(M).
\ee
Thus any sequence of indefinitely increasingly long geodesics, like the prime
geodesics found by Gromoll-Myer \cite{GlMy}, have minimizing index
approaching infinity.  This is also known to be true of the Morse
Index which will be discussed later in Section~\ref{sectenergy}.

%------
\sect{$1/2$ Geodesics and the Covering Spectrum} \label{sectcov}

%-----------------

In this section we produce a wealth of $1/2$ geodesics in
length spaces which are not simply connected.

\begin{lemma} \label{w1} 
A closed geodesic which is the shortest among all noncontractible closed
geodesics is a $1/2$ geodesic.
\end{lemma}

This lemma is a consequence of the following one applied to
the universal cover.

\begin{lemma}\label{coverhalf}
If $\tilde{M}$ is a covering space of $M$ and $c$
is the shortest curve which lifts nontrivially to
$\tilde{M}$, then $c$ is a $1/2$ geodesic. 
\end{lemma}

\Pf
If $c:S^1\to M$ is not
a $1/2$ geodesic, then $\exists t_0>0$ such that
$
d(c(t_0), c(t_0+\pi) )< L/2.
$
So we can join $c(t_0)$ to $c(t_0+\pi/L)$ by a geodesic
segment $\eta$ of length $<L/2$. Let  $c_1$ to be the
curve created by taking $c$ restricted to $[t_0, t_0+\pi]$
followed by $\eta^{-1}$ and $c_2$ to be $\eta$ followed by
$c$ from $t_0+\pi$ wrappping around to $t_0$.  Since both
$c_i$ are shorter than $c$, they must lift trivially to $\tilde{M}$.
This forces $c$ to lift trivially as well, contradicting
the hypothesis.
\qed

Note that closed geodesics which are minimizers in their homotopy
classes are not necessarily $1/2$ geodesics.  This can be seen
by looking at the geodesics in Figure~\ref{tosegment} or by
considering iterated geodesics in a torus.  The covering
spaces which unwrap these geodesics, also unwrap shorter geodesics,
thus these geodesics do not satify the hypothesis of Lemma~\ref{coverhalf}.

Geodesics which do satisy the hypothesis of Lemma~\ref{coverhalf}
were studied in \cite{SoWei}.  Their lengths correspond to the
elements of the covering spectrum, which is defined 
using a special selection of covering spaces [Defn 3.1, Theorem 4.12 
\cite{SoWei}].  We can now improve this theorem.

\begin{theorem}\label{covspec}
If $X$ is a compact space with a simply connected
universal cover, then 
\be
2CovSpec(X) \subset L_{1/2}(X).
\ee
\end{theorem}

\Pf
Theorem 4.12 of \cite{SoWei} ,
stated that $2CovSpec(X)\subset L(X)$.
A key step in the proof is Lemma 4.9 of \cite{SoWei}, where 
one takes any element $\delta \in CovSpec(X)$ and produces a
a corresponding curve $c$ of length $2\delta$ which 
satisfies the hypothesis of Lemma~\ref{coverhalf} above.
\qed

%------------------
\sect{Systoles and $1/2$ Geodesics} \label{secthalf}
%-------------------

Recall that the systole of a manifold, $sys(M)$, is the length of the
shortest noncontractible closed geodesic (c.f.\cite{CrKt}).  This definition
easily extends to any compact length space with a universal
cover or a positive injectivity radius.  Without
a positive injectivity radius the systole may be $0$ (see   
Example~\ref{Hawaii}).

Combining
Lemma~\ref{coverhalf} applied to the universal
cover of $M$ we immediately
obtain the following lemma:

\begin{lemma} \label{systole}
If $M$ is a compact length space with a universal cover, then
the shortest $1/2$ geodesic has length 
\be
sys(M) \in L_{1/2}(M)\subset (0, 2\, diam(M)]
\ee
thus $\min L_{1/2}(M) \le sys(M)$.
\end{lemma}

Note that it is quite possible that this is a strict inequality 
as the shortest $1/2$ geodesic could be contractible and wrapped
around some small ``knob'' (c.f. \cite{CoHng}).

There is a significant body of research providing upper bounds on the systole
of various surfaces, and thus also $\min L_{1/2}$. 
See for example,  Croke and Katz's recent 
survey article \cite{CrKt}.
Croke and Katz combine an inequality of Gromov \cite{G2} with
a theorem by Pu \cite{Pu}, to obtain the following proposition
which we rephrase using Corollary~\ref{systole}.

\begin{prop} [Pu] \label{Pu}  
If $M^2$ is not diffeomorphic to a sphere, then
\be
(\min L_{1/2}(M))^2 \le sys(M)^2 \le \pi Vol(M)/2
\ee 
and when
equality holds, 
$M^2$ is the standard $\RR P^2$ with constant sectional curvature.
\end{prop}

We may also rephrase Loewner's result (c.f. \cite{CrKt}).

\begin{prop} [Loewner] \label{Low1} 
If $M^2$ is diffeomorphic to a torus, then
\be
(\min L_{1/2}(M))^2 \le sys(M)^2 \le 2 Vol(M)/\sqrt{3}
\ee 
and when
equality holds, $M$ is a skewed flat torus with a $120$ degree angle. 
\end{prop}

In the following example the equality is only achieved on a singular
manifold, so it is necessary to use our compact length
definition of a $1/2$ geodesic. \cite{Bav}\cite{Sak}

\begin{prop} [Bavard, Sakai] \label{BavSak} 
If $M^2$ is diffeomorphic to $\RR P^2 \# \RR P^2$, then
\be
(\min L_{1/2}(M))^2 \le sys(M)^2 \le \pi Vol(M)/2^{3/2}.
\ee 
In this case equality doesn't hold on a manifold, but
rather on a singular space formed by gluing together two
moebius strips of constant sectional curvature $1$,
with width $\pi/2$ and central curve of length $\pi$ along
a singular circle.  The singular circle is a geodesic
in the metric space sense and achieves the minimal length,
$\sqrt{2}\pi=\pi Vol(M)/2^{3/2}$. 
\end{prop}

%Sakai, Takashi(J-OKAY)
%A proof of the isosystolic inequality for the Klein bottle.
%Proc. Amer. Math. Soc. 104 (1988), no. 2, 589--590.
%
%Baravrd refer in this ams review

In the case of manifolds diffeomorphic to $S^2$, all geodesics
are contractible and the situation is much more complicated.
Calabi and Croke have conjectured that on a surface
diffeomorphic to a sphere,
\be
\min L \le \sqrt{12} \, \sqrt{Vol(M))},
\ee
and an example achieving this inequality appears in \cite{Cr}.

%Croke, Christopher B.(1-PA) 
%Area and the length of the shortest closed geodesic. 
%J. Differential Geom. 27 (1988), no. 1, 1--21.

The strongest result in this direction is by Rotman \cite{Ro}.
%''The length of a shortest closed geodesic and the 
% area of a $2$-dimensional sphere'', to appear in Proceedings of the AMS.

\begin{prop} [Rotman] \label{Ro}
If $M^2$ is diffeomorphic to a sphere then
\be
\min L \le 4\sqrt{2}\sqrt{Vol(M)}.
\ee
\end{prop}

The best estimate based on diameter is in \cite{NaRo} and 
independantly \cite{Sab}:

\begin{prop} [Nabutovsky-Rotman, Sabourou] \label{rotnab}
If $M^2$ is diffeomorphic to a sphere then
\be
\min L \le 4 \, Diam(M).
\ee
\end{prop}

Note that neither Proposition~\ref{Ro} nor~\ref{rotnab}
provide bounds on $\min L_{1/2} \le \min L$.
It would be interesting to examine their proofs and see whether
their techniques would provide upper bounds on the
various $\min L_{1/k}$ [Problem~\ref{rotnabq}].

\begin{rmk} \label{volsys}
Note that if we were to try to extend these volume estimates to
compact length spaces then we would need a measure and a dimension
for the spaces.  One might study compact length spaces with finite
second Hausdorff measure.  See Problem~\ref{volsysq}.
\end{rmk}

There are many beautiful results estimating $\min L(M)$ for manifolds
with curvature bounds, but discussion of such results and their 
relation to the $L_{1/k}$ must be posponed to future papers.

%ADD GROMOV THM???  

%---------------------
\sect{Estimating the Length of the Shortest Closed Geodesic} \label{sectmin}
%-------------------

In this section we discuss how $1/k$ geodesics may be used
to estimate the length of the shortest closed geodesic in a 
Riemannian manifold or compact length space.  

Gromov \cite{G2}
has conjectured that for a compact Riemannian manifold $M$,
\be
\min L(M) \le c(n)Vol(M^n)^{\frac{1}{n}}
\ee
and another well-known conjecture is that
\be
\min L(M) \le c(n)Diam(M^n).
\ee
In fact Rotman suggests that
\be \label{rotmandiam}
\min L(M) \le 2 Diam(M^n)
\ee
and there are no known counterexamples.
Note that (\ref{rotmandiam}) is trivially true
when the manifold is not simply connected (c.f. covspec lemma
and Remark~\ref{L1/2nonempty}).

\begin{remark} \label{L1/2nonempty}
It follows from Lemma~\ref{diam}, that if $L_{1/2}(M)$ is nonempty,
then 
\be
\min \, L(M) \le \min \, L_{1/2}(M) \le 2diam(M).
\ee 
However, the author has recently been informed that Wing Kai Ho 
has produced examples of smooth manifolds diffeomorphic to
$S^2$  which have no $1/2$ geodesics  \cite{Ho}.
In Problem~\ref{whennonempty},
we ask what properties can be imposed on a manifold that would
guarantee the existence of a $1/2$ geodesic.
\end{remark}

To provide an estimate on $\min L$ it is necessary to define the following
quantity:

\begin{defn} \label{minindM}
The minimizing index, $minind(M)$, of a compact length space, $M$, 
is the smallest $k$ such that there is a geodesic of 
minimizing index $k$.  
\end{defn}

\begin{theorem} \label{minindMd}
If $M$ is a compact length space then
\be
\min \, L(M) \le minind(M) diam(M).
\ee
\end{theorem}

\Pf
Setting $k=minind(M)$, we know
We know $\min \, L(M) \le  \min \, L_{1/k}(M) \le k diam(M)$
by Lemma~\ref{diam}
\qed

The following theorem might help one find counterexamples
to overly sharp conjectures regarding $\min L(M)$.
See Problems~\ref{minindMdq1} and~\ref{minindMdq2}.

\begin{theorem} \label{minindMinj}
If $M$ is a compact length space and $k=minind(M)$ then 
\be
\min \, L(M) \ge \min \{ k \, injrad(M), \min \, L_{1/k}(M)\}.
\ee
\end{theorem}

\Pf
By Lemma~\ref{injLk},
if $M$ has $injrad(M)\ge i_0>0$ and minimality
index $k$, then taking $L_0=i_ok$ we have
\be
L(M)\cap [0,L_0]= L_{1/k}(M) \cap [0, L_0].
\ee
Thus $\min \, L(M)$ either $=\min \, L_{1/k}(M)$
or it is $> i_ok$.  
\qed

When $M$ is a manifold we can use Klingenberg's old
lemma to estimate its minimizing index as follows
\cite{K1} (c.f. \cite{doC}):

\begin{lemma} [Klingenberg] \label{klinglem}
Let $M$ be a compact Riemannian manifold.
If $x$ is the closest cut point to $y$, then either
$x$ is a conjugate point of $y$, or there are exactly
two geodesics from $y$ to $x$ and they meet at $180^\circ$.
If $x$ and $y$ are cut points such that $d_M(x,y)=injrad(M)$,
then either they are conjugate points or there are exactly
two geodesics from $x$ to $y$ and together they form a closed
geodesic.
\end{lemma}

\begin{coro} \label{klingcor}
If $M$ is a compact Riemannian manifold with no conjugate points
then the shortest closed geodesic is a $1/2$ geodesic of
length $2 injrad(M)$.  So $minind(M)=2$ and $min L= 2 injrad(M)$.
\end{coro}

In Problem~\ref{conjpt} we ask for an appropriate extension of
the definition of conjugate point to compact length spaces
which might allow one to extend Corollary~\ref{klingcor}.

Note that Klingenberg applied his lemma to manifolds with
negative sectional curvature as these spaces have no conjugate
points.  We suggest in Problem~\ref{klingcat} that Corollary~\ref{klingcor}
might extend to CAT(0) spaces.

In general, however, the author leaves the discussion of the
rich literature concerning the length spectrum and sectional
curvature bounds out of this paper.  We will discuss Ricci curvature
bounds as such bounds lead to applications relating to Gromov-Hausdorff
convergence.  

\sect{Convergence without Sudden Disappearances} \label{sectnodis}

In this section we prove our main convergence theorem
and present some simple illustrative examples.

\begin{theorem} \label{nodis}    \label{conv} 
If $M_i \to M$ in the GH sense then $L_{1/k}(M_i)$
converges to a subset of $L_{1/k}(M)\cup \{0\}$ in the Hausdorff sense.
That is, for all $\epsilon, R > 0$, there exists $N_\epsilon \in \NN$
such that 
\be
L_{1/k}(M_i)\cap [0,R] \subset 
T_\epsilon( L_{1/k}(M) \cup \{0\} ) \qquad \forall i\ge N_\epsilon.
\ee
\end{theorem}

Recall that in Figure~\ref{boobah} we gave an example
with suddenly appearing
$1/2$ geodesics, which are not the limits of such 
a sequence of $L_i$ or even
$L_i$ selected from $L(M_i)$.  

\Pf
Suppose, on the contrary, 
there exists $\epsilon$ such that for a subsequence of the $i$
there are 
\be
L_i \in L_{1/k}(M_i) \setminus 
T_\epsilon( L_{1/k}(M) \cup \{0\} ).
\ee
By Lemma~\ref{diam} and the fact that Gromov-Hausdorff convergence
implies that $diam(M_i)\to D=diam(M)$,
$L_i \in [0,kdiam(M_i)]\subset[0, 2kD]$. 
So a further subsequence must converge, $L_i\to L_0$, where
\be \label{nodis2}
L_0 \in (0,2kD] \setminus L_{1/k}(M).
\ee  
Thus there exists corresponding $1/k$ geodesics 
$\gamma_i: S^1 \to M_i$ of speed $L_i/(2\pi)$.  
By Grove Petersen's Arzela-Ascoli
Theorem \cite{GrPet}, 
 we know that a subsequence of the $\gamma_i$ converges
to a curve $c:S^1 \to M$ of length $L_0$.  Furthermore
\begin{eqnarray}
d_M(c(t-pi/k, c(t+\pi/k))  
&= & \lim_{i \to \infty} d_{M_i}(\gamma_i(t-\pi/k), \gamma_i(t+\pi/k)) \\ 
&= & \lim_{i \to \infty} L_i/k  = L_0/k,  
\end{eqnarray}
so $c$ is either a $1/k$ geodesic or it is trivial.
This contradicts (\ref{nodis2})
\qed

\begin{rmk}
The same proof could be used to show that if $\gamma_i$ all have
$injrad(\gamma_i)\ge r_0$, then their limit does as well.
\end{rmk}

The following two corollaries are immediate:

\begin{coro} \label{minindtoinfty} 
Suppose $L_i \in L(M_i)$ and $L_i\to L_\infty \notin L(M)$ where
$M$ is the Gromov-Hausdorff limit of the $M_i$.  Then 
the geodesics $\gamma_i:S^1\to M_i$ 
of length $L(\gamma_i)=L_i$ have $minind(\gamma_i)$ diverging to infinity.
\end{coro}

\begin{coro} \label{L1/kcompact} 
Given any compact length space $M$, $L_{1/k}(M)\cup \{0\}$ is compact.
\end{coro}

%Note that $L(M)\cup \{0\}$ is not necessarily compact [Example~\ref{notcomp}].

\begin{example} \label{torus1} 
Recall $M_j=S^1_{\pi}\times S^1_{\pi/j}$, the flat tori of 
Example~\ref{torusH} that converged to
a circle whose length spectra had disappearing lengths.  
In Example~\ref{torus0}, we showed $L_{1/(2k)}(M_j)$ is  
the union:
\be
\{\sqrt{(2\pi a)^2 
+ (2\pi b/j)^2 \, }: a,b = 1,2,...k \}
\cup
\{2\pi b/j: b = 1,2,...k \}
\cup
\{2\pi a: a = 1,2,...k \}
\ee
and $L_{1/(2k-1)}(M_j)=L_{1/(2k)}(M_j)$.
As $j$ diverges to infinity, $L_{1/(2k)}(M_j)$, converges in
the Hausdorff sense to the union
\be
\{\sqrt{(2\pi a)^2 + (0b)^2 \, }: a,b = 1,2,...k \}
\cup
\{ (0b: b = 1,2,...k \}
\cup
\{2\pi a: a = 1,2,...k \}
\ee
which is $=\{2\pi a: a=0,1,2,...k\}=L_{1/(2k)}(S^1)$.  
\end{example}

%\begin{note}  REMOVE THIS NOTE AND POSPONE UNTIL WORK WITH ROTMAN
%The shortest length in the length spectrum of a manifold need
%ot have minimizing index $=1/2$ even when $L_{1/2}$ is nonempty.
%This might be seen by taking a standard sphere and instead
%of putting two handles, use a bulging ridge of length $\pi/2$
%that holds a geodesic of length $\pi$.   $L_{1/2}$ will still 
%include $2\pi$ but a closed geodesic $\gamma_i$ of length $L_i$ 
%approximately $\pi$ will run around the ridge.  If we take a limit 
%as the knobs get smaller and smaller,
%this converges in the GH sense to the standard $S^2$, but
%these $L_i \to \pi \notin L(S^2)$.  Thus the minimizing index of these
%$L_i$ is going to infinity.  Note that I use a ridge instead of 
%handles or knobs because handles and knobs have shorter closed geodesics
%running around them that are in $L_{1/2}$ and converge to $0$.
%However, it is not clear how to show that the geodesic running around
%the ridge is the shortest, so this is just a note and not an example.
%\end{note}

\begin{example} \label{ellipsoid1} {\em
Recall the sequence of ellipsoids, $M_j^2$, from Example~\ref{DD}
converging to a doubled disk, $Y$.

Note that the curves $h_j(t)=(cos(t),sin(t),0)$ mapped into $M_j$
are closed geodesics.  Their pointwise limit as $c_j \to 0$
is $h_\infty(t)=(cos(t), sin(t), 0)$ mapped into $M_\infty$.  
Note that $h_\infty$ is parametrized by arclength but 
\be
d_{M_\infty}(h_\infty(t-\epsilon), h_\infty(t+\epsilon))=
2sin(\epsilon)<2\epsilon.
\ee
So $h_\infty$ is not a closed geodesic.
{\em Thus there exists no uniform lower bound on 
the minimizing index of the $h_j:S^1\to M_j$.}}
This can also be seen using the recent work of \cite{ItKi}
\end{example}

The next example demonstrates why it is necessary to study
$1/k$ geodesics rather than just smooth regular polygons
that are minimizing between only k specific regularly spaced points
instead of any collection of k regularly spaced points.

\begin{example} \label{hexagon1}
In Figure~\ref{squarefig} we see
$M_\epsilon \subset \EE^3$, the boundary of the $\epsilon$
tubular neighborhood around a 
flat solid regular square $Z \in \EE^2\times\{0\}$.  
For $\epsilon$ sufficiently
small we can see that the geodesic, $\gamma_\epsilon:S^1\to M_\epsilon$ 
running around the equator looks
almost like a square.  If one chooses a specific regularly
spaced selection of four points on $\gamma_\epsilon$ each of which
is close to the corner of the square, one sees that $\gamma_\epsilon$
is minimizing between these points.  However $\gamma_\epsilon$ is
not a $1/4$ geodesic.
%\psdraft
\begin{figure}[htbp]
\centering
\includegraphics[height=1in ]{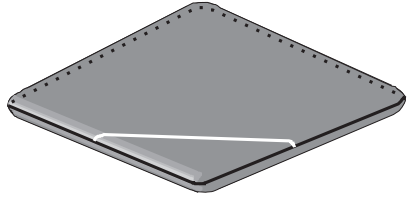}
\caption{The black closed geodesic, $\gamma_\epsilon$, is minimizing between the four corners
but not their midpoints as indicated by the white geodesic segment. }
\label{squarefig}
\end{figure}
%\psfull
%hexagon figure

As in Example~\ref{ellipsoid1}, the limit space as $\epsilon$
approaches $0$ is a doubled copy of $Z$ glued to itself along the
square boundary.  The square boundary is only 
piecewise minimizing between
the corners and is not a closed geodesic.
Thus the $\gamma_\epsilon$ do not converge to a 
closed 
geodesic in the limit
space.  
\end{example}

\begin{remark} \label{bangertrmk}
If $M_i$ converge to the standard $S^n$ smoothly, then Bangert
proved that the prime geodesics in $M_i$ either have lengths
converging to $2\pi$ or to $\infty$ \cite{Bng}.  Note that 
the prime geodesics $\gamma_i$ whose lengths diverge to infinity, 
have minimal indices also diverging to infinity by Lemma~\ref{diam}.

Bangert's Theorem does not extend to $M_i$ converging to $S^n$
in the Gromov-Hausdorff sense.
If we take $M_j=S^n\times S^1_{1/j}$
and prime geodesics wrapping once around the equator of $S^n$ while
wrapping $j$ times around the $S^1_{1/j}$.  These all have length $4\pi$. 
Also the geodesics in Example~\ref{tosegment} are prime geodesics
converging to a length $<\pi$.  
\end{remark}

In Problem~\ref{bangertq} we ask whether Bangert's Theorem
described here holds when one assumes the $M_i$ converge
in the Gromov Hausdorff sense with a uniform positive
lower bound on injectivity radius.

%==============================
\sect{Gap Theorems and Ricci Curvature} \label{almrig}
%==============================

% stability vs almost rigidity and cosmos

In this section we apply the length spectrum convergence theorem
[Theorem~\ref{nodis}] to force the existence of gaps in the
length spectrum of certain manifolds with Ricci curvature bounds.  
We begin by rephrasing Theorem~\ref{nodis} as a gap theorem:

\begin{theorem}\label{gapthmk}
Fix a compact length space, $M$, and choose any $\epsilon>0$
$k\in \NN$,
then there exists $\delta_{\epsilon, k, M}$ such that if
\be \label{gapthm1}
(a,b) \cap L{1/k}(M)=\emptyset
\ee
then
\be
[a+\epsilon, b-\epsilon]\cap L_{1/k}(N)=\emptyset
\ee
for all compact length spaces $N$ such that
$d_{GH}(N,M)<\delta_{\epsilon, k, M}$.
\end{theorem}

\Pf 
Suppose on the contrary that there exists $\epsilon>0$, $k\in \NN$ 
and $N_i$ converging to $M$ with 
\be
L_i \in [a+\epsilon, b-\epsilon]\cap L_{1/k}(N_i).
\ee
Then by Theorem~\ref{nodis}, a subsequence of the $L_i$
converge to some $L\in L_{1/k}(M)$.  Since
$L_i\in [a+\epsilon,b-\epsilon]$, 
so is $L$ which contradicts (\ref{gapthm1}).
\qed

The next gap theorem refers directly to the length spectrum.

\begin{theorem} \label{gapi}
Fix a compact length space, $M$, and choose any $\epsilon>0$
and $b>0$,
then there exists $\delta_{\epsilon, b, M}$ such that if
\be \label{gapi1}
(ai_0,bi_0) \cap L(M)=\emptyset
\ee
then
\be \label{gapi2}
[ai_0+\epsilon, bi_0-\epsilon]\cap L(N)=\emptyset
\ee
for all compact length spaces $N$ such that
$d_{GH}(N,M)<\delta_{\epsilon,b, M}$
with $injrad(N)\ge i_0$.
\end{theorem}

\Pf
First note that (\ref{gapi1}) implies that
\be
(ai_0,bi_0) \cap L_{1/k}(M)=\emptyset \qquad \forall k\in \NN.
\ee
Now we choose $k \ge b$, and apply Theorem~\ref{gapthmk}
for that $k$, and take
$\delta_{\epsilon,b, M}:=\delta_{\epsilon, k,M}$,
which implies that
 \be
[a+\epsilon, b-\epsilon]\cap L_{1/k}(N)=\emptyset.
\ee
Restricting to $N$ with $injrad(N)\ge i_0$ and applying
Lemma~\ref{injLk}, we get (\ref{gapi2}).
\qed

Notice if one happens to take a sequence of $N_i\to M$
whose injectivity radii converge to $0$, we can still
apply Theorem~\ref{gapi} but the gaps slide over towards
$0$ and shrink.  This is seen to be exact in Example~\ref{torusH}.

When one has $C^2$ convergence of the manifolds, Ehrlich has proven
the injectivity radii converge, in which case Theorem~\ref{gapi}
is significantly stronger and basically already known \cite{Eh}.
%(cf Weinstein REF).

It is crucial to understand that even with smooth convergence
we do not get uniform $\delta$ depending only on $\epsilon$.
They will always depend on the manifold itself.  Otherwise
we would never have suddenly appearing geodesics.  That is,
if $\delta_{M,b, \epsilon}$ did not depend on $M$, then
take $M_i$ converging smoothly to $Y$ as in Figure~\ref{boobah}
and $b$ twice the length of the suddenly appearing geodesic, then
for $d_{GH}(M_i, Y)<\delta_{b, \epsilon}$ we'd contradict
the existence of that suddenly appeared geodesic.  In fact
$\delta_{M_i,b, \epsilon} > (b/2) \epsilon d_{GH}(M_i, Y)$,
so that $N$ which are much closer to $M_i$ than $M_i$ is
to $Y$ will have a gap such that there is no geodesic
of length $b/2$.

%For a specific manifold $M$, it can be quite simple to estimate
%explicit values for $\delta_{M,k,\epsilon}$.  
%This will be
%discussed in \cite{So12} where we develop the concept of
%an almost $1/k$ geodesic.

Thereoms~\ref{gapthmk} and~\ref{gapi} can now be applied to
a number of stability  results to prove the gap theorems
mentioned in the introduction.

Recall that Bishop proved that any Riemannian manifold
$M^n$ with $Ricci\ge (n-1)$ and $Vol(M^n)=Vol(S^n)$ is isometric
to the sphere \cite{Bi}.  The stability theorem for this rigidity
result was proven by Colding \cite{Co1}:

%volume stability
%Colding, Tobias H.(1-NY-X) 
%Shape of manifolds with positive Ricci curvature. 
%Invent. Math. 124 (1996), no. 1-3, 175--191.

%Anderson, Michael T.(1-SUNYS) 
%Convergence and rigidity of manifolds under Ricci curvature bounds. 
%Invent. Math. 102 (1990), no. 2, 429--445.

\begin{prop} [Colding]\label{coldingvol}
If $N^n$ has $Ricci(N^n)\ge -(n-1)$, then for all $\epsilon>0$,
there exists $\delta_{\epsilon,n}>0$ such that 
\be
Vol(N^n) \ge Vol(S^n)-\delta_{\epsilon,n}
\ee
implies $d_{GH}(N^n, S^n)<\epsilon$.
\end{prop}

This proposition combined with Theorem~\ref{gapthmk}
and the length spectrum of $S^n$ in Example~\ref{sphere0}
implies the following:

\begin{prop} \label{GHSL}
For all $\epsilon>0$,
and any $k\in \NN$,
there exists $\delta_{\epsilon,k,n}>0$ such that 
\be
Vol(N^n) \ge Vol(S^n)-\delta_{\epsilon,k,n}
\ee
and $Ricci(N^n) \ge -(n-1)$ 
then for all $j\in \{0,1,2,...\}$ we have:
\be
[2j \pi +\epsilon, 2(j+1)\pi -\epsilon]\cap L_{1/(2k)}(N)=\emptyset.
\ee
\end{prop}

Combining this proposition with Lemma~\ref{GHdiam} and Lemma~\ref{diam}
we obtain Theorem~\ref{volgap} which was stated in the introduction.

\begin{rmk}\label{volgapr}
Problem~\ref{volgapq} asks for precise estimates on the
estimating function in Theorem~\ref{volgap}.  Given a precise estimate, 
one would be
able to bound the volume of a manifold $N^n$ with $Ricci \ge (n-1)$
depending on the length and minimal index
of one of its closed geodesics.
\end{rmk}

Although Colding did later prove convergence in the $C^{1,\alpha}$
topology \cite{Co3}, there are manifolds, 
$M^n_i$, satisfying
$Ricci(M^n_i) \ge (n-1)$ , $Vol(M_i)\to Vol(S^n)$ whose
injectivity radii $injrad(M_i) \to 0$ [Example~\ref{football}]
So we cannot presume to improve the
length spectrum's convergence or obtain an estimate on $\min L(N^n)$
without imposing an additional condition on the injectivity radius.

%\cite{Co3}
%C^{1,\alpha} and more also 
%Colding, Tobias H.(1-NY-X) 
%Ricci curvature and volume convergence. 
%Ann. of Math. (2) 145 (1997), no. 3, 477--501.

\begin{example}  \label{football}
We now construct smooth American footballs, $M^2_j$, with $sect \ge 1$ 
(and thus $Ricci\ge 1$) whose volume $Vol(M_j)\ge a_j^2 Vol(S^2)-\epsilon_j$
with $a_j \to 1$ and $injrad(M_j)\le r_j \to 0$.
 
Start with the standard $S^2$, remove a wedge of angle $(1-a_j)2\pi <\pi/4$,
and glue the edges to themselves to get a singular manifold, $F_j$,  of
volume $a_j^2 Vol(S^2)$.  For small $r_j>0$, take 
two points $p_j, p_j'$ both $r_j/2$ away
from a singular point and  maximally far apart.  There are two distinct
geodesics running between them of length less than $r_j$.  Let $h_j$
be the distance from these geodesics to the singular point.  Then
$h_j > r_jcos(a_j\pi/2)$, the height of the Euclidean comparison triangle.

%football figure

Now if we remove balls of radius $h_j/2$ about the two singular points
in $F_j$, we can cap off these regions smoothly with caps whose $sect \ge 1$
This gives our surfaces $M_j$ and the points $p_j$ and $p_j'$ are still
cut points in $M_j$ and so $injrad(M_j)< r_j$ and 
$Vol(M_j) \ge Vol(F_j)- \pi r_j^2$.

It is not clear how the length spectrum of these $M_j$ behave.
Are there examples of $M_j$ with $\min L(M_j)\to 0$ or a disappearing
length?  [Problems~\ref{volgapmin} and~\ref{volgapdis}].
\end{example}

Myers Theorem states that any manifold $M^n$ with
$Ricci \ge (n-1)$ has $diam(M^n)\le \pi$ because
any geodesic of length $\pi$ must have a conjugate point 
(c.f. \cite{doC} \cite{My}).
Cheng's Sphere Rigidity Theorem
states that this inequality is only achieved on a sphere \cite{Chng}. 

Cheng's Theorem doesn't have
a stability theorem like Proposition~\ref{coldingvol},
as is demonstrated by Otsu's examples \cite{Ot}.  
Otsu's five dimensional manifolds satisfy the Ricci bound
and their diameter approaches $\pi$ but they converge
in the Gromov-Hausdorff sense to a singular manifold
not a sphere.  This limit space contains only two points
which are a distance $\pi$ apart.

%Myers, S. B. Riemannian manifolds with positive mean curvature. 
%Duke Math. J. 8, (1941). 401--404.

%Otsu, Yukio(J-KYUSS) 
%On manifolds of positive Ricci curvature with large diameter. 
%Math. Z. 206 (1991), no. 2, 255--264.

%Cheng SY Eigenvalue comparison theorems and its geometric applications,
% Math Z 143 (1975) 289-297

\begin{rmk}\label{diamgap}
Cheeger-Colding have proven that a manifold
with $Ricci\ge (n-1)$ and diameter close to $\pi$ is Gromov-Hausdorff close
to a spherical suspension over a subset of the manifold \cite{ChCo}.
This is called an ``almost rigidity'' result rather than a stabilit
result because they do not prove it is close to a particular metric space,
but rather than the metric bahaves in an almost rigid manner.
In Problem~\ref{diamgapq}, we question whether one can obtain
a gap theorem based on such a result.  One of the biggest
difficulties there would be turning this spherical suspension into
a length space and not just a metric space.
Then naturally one would need to know if there are any uniform
properties of the length spectrum on spherical suspensions
[Problem~\ref{sphsusp}].
\end{rmk}

%Cheeger, Jeff; Colding, Tobias H. 
% Lower bounds on Ricci curvature and the stability of warped products. 
% Ann. of Math. (2) 144 (1996), no. 1, 189--237. 

To avoid the issue arising in Otsu's examples, Colding instead
examined the radius:

\begin{defn} \label{radiusdef}
The radius of a compact metric space, $M$, is the smallest $r>0$ such 
that $M\subset \bar{B}_p(r)$ for some $p$.  In fact
\be
rad(M)=\inf_{p\in M} \sup_{q\in M} d(p,q)\le diam(M).
\ee
\end{defn}

When a manifold with $Ricci\ge (n-1)$ has radius close to $\pi$,
then every point in the manifold has a point almost maximally
distant from it, thus it is approaching the inequality in
Myer's Theorem along every geodesic in the manifold.
Colding proved the following stability result \cite{Co2}

%volume continuity and radius stability
%Colding, Tobias H.(1-NY-X) 
%Large manifolds with positive Ricci curvature. 
%Invent. Math. 124 (1996), no. 1-3, 193--214.

\begin{prop} \label{radius} [Colding]
Given $n\ge 2$ and $\epsilon>0$ there exists $\delta(n, \epsilon)>0$
such that if $N^n$ is a compact Riemannian manifold whose
$Ricci(M) \ge (n-1)$ and $rad(M)> \pi -\delta$ then 
$Vol(M)>Vol(S^n)-\epsilon$.
\end{prop}

Combining Proposition~\ref{radius} with Theorem~\ref{gapthmk} we obtain:

\begin{theorem}\label{radgap} 
There exists a function $\Psi:\RR^+\times \NN \times \NN \to \RR^+$ such
that $\lim_{\delta \to 0}\Psi(\delta, k,n)=0$ such that
if $N^n$ is a compact Riemannian manifold with
\be
rad(N^n) \ge \pi-\delta
\ee
and $Ricci(N^n) \ge (n-1)$ then
\be
L_{1/(2k)}(M^n) \subset [0,\epsilon_k)\cup (2\pi-\epsilon_k, 2\pi+\epsilon_k) 
\cup \cdots (2k\pi-\epsilon_k, 2k\pi+\epsilon_k)
\ee
for $\epsilon_k=\Psi(\delta,k,n)$.
\end{theorem}

\begin{rmk} \label{footballrad}
Note that Example~\ref{football} also has $rad(M_i)\to rad(S^2)$, so
here we also have no lower bound on injectivity radius and cannot
directly conclude a stronger convergence of the length spectrum.
See Problems~\ref{radgapmin} and~\ref{radgapdis}.
\end{rmk}

Gromov proved that any Riemannian manifold
$M^n$ with $Ricci\ge 0$ and first Betti number satisfying $b_1(M)=n$ 
is isometric to a torus \cite{G1}.
The corresponding stability theorem  is hidden
in Colding's proof
that any $M^n$ with $b_1(M^n)=n$ and $Ricci\ge -(n-1)\epsilon $
is homeomorphic to $\TT^n$ if $\epsilon$ is sufficiently small
\cite{Co3}.

\begin{prop} [Colding]
For any $\epsilon>0$ there exists $\delta_\epsilon>0$ such that if
$M^n$ has $b_1(M^n)=n$ and $Ricci \ge -(n-1)\delta_\epsilon$ then
\be
d_{GH}(M^n, \TT^n) < \epsilon.
\ee
\end{prop}

Combining this with Theorem~\ref{gapthmk} and the fact that
\begin{eqnarray} \label{GHTL1}
L_{1/(2k)}(\TT^n) &=&\{L^k_1, L^k_2,...L^k_{m(k)}\} \\
& = & \{\sqrt{\, \sum_{m=1}^n a_n^2 \,}: a_m=0,1,2,...k\} 
\setminus \{0\}
\end{eqnarray}
we get the following:

\begin{theorem} \label{torusgap}
For all $\epsilon>0$,
and any $k\in \NN$,
there exists $\delta_{\epsilon,k}>0$ such that 
if $N^n$ is a Riemannian manifold with
\be
Ricci(N^n) \ge -(n-1)\delta_{\epsilon,k}
\ee
and $b_1(N^n)=n$ then
then for all $j\in \{0,1,2,...m(k)\}$ we have:
\be
[L^k_j +\epsilon, L^k_{j+1} -\epsilon]\cap L_{1/(2k)}(N)=\emptyset,
\ee
where $0=L^k_0 < L^k_1<L^k_2<...$ is given in (\ref{GHTL1}).
\end{theorem}

In Problems~\ref{torusgapmin} and~\ref{torusgapdis} we ask for
relevant examples with arbitrarily small $\min L$ or
disappearing lengths.  Examples with injectivity radius approaching
$0$ have been described by Colding.  One begins by removing balls of 
radius $1/2$ from a flat torus, gluing in very flat cones, and then smoothly
capping them off carefully to keep the injectivity radius exactly
as in Example~\ref{football} and finally smoothing the boundaries of
the balls which adds the slightly negative curvature.  

%In fact Colding also proved that when $N^n$ is sufficiently close
%to $\TT^n$ it is homeomorphic to $\TT^n$, in which case we know
%the fundamental group is nontrivial and 
%\be
%\min L(M) \le \min L_{1/2}(M).
%\ee
%WITHOUT AN INJECTIVITY RADIUS BOUND THIS IS USELESS.

\begin{rmk}\label{cosmos}
In \cite{So1}, the author has proven another stability theorem,
that a locally almost isotopic manifold with $Ricci \ge -(n-1)H$ is 
Gromov-Hausdorff close to an isotopic manifold.  When the manifold
is compact it close to a Riemannian manifold homothetic to a sphere.
Thus we would also get a gap theorem for such a manifold.  Since the
definition of locally almost isotopic is complicated, we do not 
give the complete explanation here.  See Problem~\ref{cosmosq}
\end{rmk}

%--------------------------------
\sect{Openly $1/k$ Geodesics}  \label{sectopen}
%---------------------------------------

In this section we define a new collection of geodesics and
lengths which behaves a bit better than $1/k$ geodesics on
manifolds.

\begin{defn} \label{defopen}
An openly $1/k$ geodesic $\gamma:S^1\to M$ is a $1/k$ geodesic
which has $injrad(\gamma)>L(\gamma)/k$.
\end{defn}

The great advantage of an openly $1/k$ geodesic is the following
lemma.

\begin{lemma} \label{oldunique}
If $M$ is a compact Riemannian manifold and $\gamma$
is an openly $1/k$ geodesic, then it is uniquely determined
by any collection of $k$ evenly spaced points
up to reparametrization by an isometry of $S^1$.
\end{lemma}

\Pf This follows from the fact that if $\gamma$ is minimizing
on $[a,b]$ then it is uniquely determined on $[a,c]$ for any
$c\in (a,b)$ (c.f. \cite{doC}).
\qed

Lemma~\ref{oldunique} does not hold on a compact length space:

\begin{example}\label{doubledsquare}
Let $Y$ be a graph with four
ordered vertices, $\{v_1, v_2, v_3, v_4=v_0\}$, and two unit edges 
$e^+_i$ and $e^-_i$ between $v_i$ and $v_{i+1}$ for $i=0,1,2,3$.  
Let $\gamma$ be the  geodesic which traverses $e^+_1$, $e^-_2$, $e^+_3$,
and $e^-_4$.  It is in fact a $1/2$ geodesic and thus an openly
$1/4$ geodesic.  However, it is not uniquely determined by
the 4 evenly spaced points $v_1, v_2, v_3, v_4$ as there is
another geodesic sharing those points which traverses
 $e^-_1$, $e^+_2$, $e^-_3$,
and $e^+_4$. 
\end{example}

Lemma~\ref{oldunique} does not hold if one only assumes
$\gamma$ is a $1/k$ geodesic instead of openly $1/k$
as we see in the next example:

\begin{example} \label{chiclet}
If we take a Riemannian manifold depicted in Figure~\ref{chicletfig},
\be
M_\epsilon=\partial T_\epsilon([-1,1]\times[-1,1]\times\{0\})\subset \EE^3,
\ee
and choose the points 
\be
p_j=((1+\epsilon)cos(j\pi/2), (1+\epsilon)sin(j\pi/2), 0)
\ee
then we claim the piecewise geodesic, $\gamma$, which runs minimally
with positive $z$ from $p_0$ to $p_1$, negative $z$ from
$p_1$ to $p_2$, positive $z$ from $p_2$ to $p_3$ and negative $z$
from $p_3$ to $p_0$, is a $1/4$ geodesic.  

%\psdraft
\begin{figure}[htbp]
\centering
\includegraphics[width=6in ]{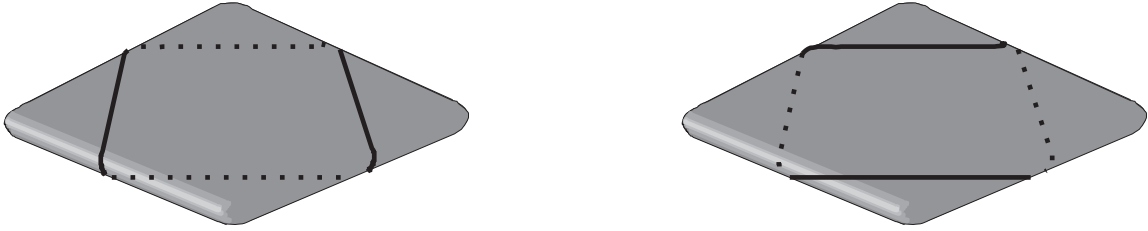}
\caption{Here we have two copies of $M_\epsilon$.  The geodesic
$\gamma$ is depicted on the right.}
\label{chicletfig}
\end{figure}
%\psfull

To prove we show that in fact $\gamma$ is actually minimizing
on $[j\pi/2 -s, j\pi/2+s]$ for any $s<\pi/2$.  First note 
the $z$ components of $\gamma(j\pi/2-s)$ is the
negative of $\gamma(j\pi/2+s)$.  So if $\sigma$ runs minimally
between these two points it must have an $s_0$ where its $z$
component is $0$.  By symmetry $\sigma(s_0)$ must be located
at $p_j$, thus $\sigma$ must agree with $\gamma$.
So $\gamma$ is actually minimizing 
between any $t$ and $t+2\pi /4$!  %add figure square

Note that $\gamma$ is also not uniquely determined by the $p_i$ because
there is another geodesic which is the reflection of $\gamma$
through the $xy$ plane running through the same four points
also depicted on the left in Figure~\ref{chicletfig}.  
\end{example}

We now develop the theory of openly $1/k$ geodesics.

\begin{defn}
Let $L^{open}_{1/k}(M)$ be the collection of lengths of
openly $1/k$ geodesics.
\end{defn}

The following lemma is an easy exercise:

\begin{lemma} \label{back}
For any $k>2$ we have
\be
L_{1/(k-1)}(M) \subset L^{open}_{1/k} \subset L_{1/k}(M)
\ee
and $L^{open}_{1/2}(M)=\emptyset$. 
\end{lemma}

In particular there are no openly $1/2$ geodesics.
Lemma~\ref{back} combined with Theorem~\ref{bigcup} immediately implies:

\begin{theorem}~\label{openbigcup}
On any compact length space, $M$,
\be
L(M)=\bigcup_{k=3}^\infty L^{open}_{1/k}(M).
\ee
\end{theorem}

\begin{defn}
Let the openly minimizing index, $opind(\gamma)$, 
of a geodesic, $\gamma$, be the
smallest $k$ such that $\gamma$ is an openly $1/k$ geodesic.

Let $opind(M)=\min\{opind(\gamma): \gamma:S^1\to M\}$.
\end{defn}

Lemma~\ref{back} immediately implies:

\begin{lemma}\label{minop}
\be
minind(\gamma) \le opind(\gamma)\le minind(\gamma)+1.
\ee
\end{lemma}

Note that in the flat torus and in the sphere
all closed geodesics have $minind(\gamma)<opind(\gamma)$
[Examples~\ref{torus0} and~\ref{sphere0}].
Manifolds with this property are of significant
interest because we are able to bound the open index
in  Theorems~\ref{openenergy} and~\ref{morsethm} below.  
See Problem~\ref{minopind}.

Theorem~\ref{nodis} and Lemma~\ref{back} together imply:

\begin{theorem}\label{opnodis}
If $M_j \to M$ in the GH sense then 
\be
\lim_{j\to\infty} L^{open}_{1/k}(M_j) \subset L_{1/k}(M).
\ee
\end{theorem}

We do cannot improve this to 
\be
\lim_{j\to\infty} L^{open}_{1/k}(M_j) \subset L^{open}_{1/k}(M).
\ee
as can be seen in the
following example.

\begin{example}\label{ellipsoidopen}
In Example~\ref{ellipsoid1} we demonstrated that
the equators, $\gamma_c$, of flattening ellipsoids 
$(x^2+y^2+(z/c)^2=1$ had $injrad(\gamma_c)$ varying
continuously with $c$ and converging to $0$ as $c\to 0$.

When $c=1$,
$minind(\gamma)=1/2$ and $injrad(\gamma)=\pi$.
As $c$ decreases, the injectivity radius decreases
continuously \cite{ItKi},
and at some $c_0>0$ the injectivity radius hits
$2\pi/3$ for the first time.  So for all
$c>c_0$, $\gamma_c$ is an openly $1/3$ geodesic
but $\gamma_{c_0}$ is not.

So if $c_j$ decrease to $c_0$,
\be
\min_{j\to\infty}L^{open}_{1/3}(M_{c_j}) \notin L^{open}_{1/3}(M_{c_0})
\ee
even though $M_{c_j}$ converges to $M_{c_0}$ in the $C^\infty$ and
Gromov-Hausdorff sense.
\end{example}

It is of some interest to understand what is special about
$1/k$ geodesics which are not openly $1/k$ geodesics.
On such a geodesic, $\gamma$, there is a pair of cut points which are
a distance $L(\gamma)/k$  apart.

\begin{prop} \label{klingit}
If $M$ is a compact Riemannian manifold with no conjugate points
and $\gamma$ is a $1/k$ geodesic of length $k injrad(M)$ 
which is not
an openly $1/k$ geodesic then either $k=2$ or $\gamma$ is
the iterate of a $1/2$ geodesic and $k$ is even,
\end{prop}

\Pf
If $\gamma:S^1\to M$ is a $1/k$ geodesic which is not
an openly $1/k$ geodesic, then it has a pair of cut points on it
which are a distance $L(\gamma)/k$ apart.  If $L(\gamma)/k=injrad(M)$
and $M$ has no conjugate points, then by 
Klingenberg's Lemma~\ref{klinglem},
these two points are joined by exactly two geodesics which
close up smoothly, thus either $k=2$ or $\gamma$ is an
iterated geodesic $\gamma(t)=\gamma_0(kt/2)$ with $k$ even.  
Furthermore $\gamma_0$ is a $1/2$ geodesic because 
$d(\gamma_0(t), \gamma_0(t+\pi))\ge injrad(M)=L(\gamma_0)/2$.
\qed

Proposition~\ref{klingit} is not true on metric spaces.

\begin{example} \label{Klingitex}
Let $M$ be the metric space which is a graph with two
vertices and three unit length edges each running from one vertex 
to the other.  Then $injrad(M)=1$ and any path which runs back and forth
between the endpoints with constant speed and never traverses
back on the edge it just crossed over is a geodesic.
Thus for any $k\in\NN$ $M$ has many prime $1/(2k)$ geodesics of length 
$2k injrad(M)$.
\end{example}

%-------------
\sect{Energy and Openly $1/k$ Geodesics}  \label{sectenergy}
%---------------

Here we introduce an energy  method which may be used to prove the 
existence of a $1/k$ geodesic on a given space with certain
properties, thus allowing one to estimate $minind(M)$ and
thus $\min L(M)$ via Theorem~\ref{minindMd} and Lemma~\ref{minop}.

In this section we limit ourselves to convex compact Riemannian manifolds
with boundary so that we can discuss the derivative
of a geodesic.  The convexity assumption guarantees the geodesics
won't touch the boundary.  Background material may be found in \cite{BTZ}
and \cite{Mil}.

Smoothly closed geodesics are the critical points of the energy function
on the loop space of $M$:
\be
E(c) =\int_0^1 g(c'(t), c'(t)) \,dt
\ee
It is easy to see that when we have a critical point of this energy,
one gets a smoothly closed geodesic.  Furthermore if $c$ is a smoothly
closed geodesic and is known to be minimizing on subintervals $[t_i, t_{i+1}]$
then the energy satisfies:
\be
E(\gamma) = \sum_{i=1}^N d(\gamma(t_{i+1}), \gamma(t_i) )^2/ (t_{i+1}-t_i).
\ee
So if $\gamma$ is a $1/k$ geodesic, then 
\be
E(\gamma) = \sum_{i=0}^{k-1} d(\gamma((i+1)/k), \gamma(i/k) )^2/ (1/k).
\ee

In Morse Theory one uses a uniform lower bound on injectivity
radius and makes a finite dimensional approximation of the loop space.
That is any smoothly closed geodesic of length $\le L$ can be viewed
as a critical point in 
\be
\Omega_k(M) \subset
(M)^k=M\times M \times \cdots \times M
\ee
where
\be
\Omega_k(M)=\{(x_1,..x_k): d(x_i, x_{i+1}) \le i_0 \}.
\ee
of the
energy function
\be \label{finiteenergy}
E(x_1,...x_k)=\sum_{i=0}^k d(x_i, x_{i+1})^2/(1/k)
\ee
where $k\ge L/i_0$.  Once one finds the $x_i$ which give a critical
value, you join them by the unique geodesic segments between them
to get a loop and prove that this loop is a smoothly closed geodesic.

In particular one has the following old theorem:

\begin{theorem}  \label{oldenergy} [c.f. \cite{Mil}]
If $M$ is a manifold,
given a set of length segements $r_i\in \RR^+$ we can define
\be \label{oldenergy1}
E_{\{r_1, r_2,...r_k\}}(x_1,...x_k)=\sum_{i=1}^k d(x_i, x_{i+1})^2/r_i
\ee
where $x_{k+1}=x_1$.
Then $(y_1,...y_k)$ is a smooth critical point of $E:(M)^k \to \RR$ 
iff for all $i\in \{1, 2, ...k\}$ we have:
\be \label{oldenergy2}
d(x_i, x_{i+1})/r_i=d(x_{i-1}, x_i)/r_{i-1} 
\ee
and
\be
\nabla \rho_{x_{i+1}} = -\nabla \rho_{x_{i-1}} \textrm{ at } x_i,
\ee
where $\rho_{x}(y)=d(x,y)$.
\end{theorem}

Note that $\nabla \rho_x$ is not defined at cut points of $x$.
Here however, we avoided this issue by explicitly stating that we are
at a smooth critical point.

In particular, if we study $E=E_{1/k, 1/k, 1/k,...1/k}$ on
$\Omega_k$, it is a smooth function when it's values are 
less than $L$.  So all of its critical points below $k^2L^2$,
are smooth geodesics which are minimizing between $k$ evenly spaced
points.

\begin{example} \label{hexagon2} 
Let $M_\epsilon=\partial T_\epsilon(Y)\subset \EE^3$ where
$Y$ is a flat solid regular square in $\EE^2\times \{0\}$
as in Example~\ref{hexagon1}.  For $\epsilon$ sufficiently
small we can see that the geodesic running around the equator looks
almost like a square and is the critical point of the energy in
(\ref{finiteenergy}) for $k=4$ when the $x_i$ are near the vertices
of the square.  However, it is not a minimizing geodesic between the
midpoints of the sides, and so it is not a $1/4$ geodesic.
\end{example}

Nevertheless we would like to use Theorem~\ref{oldenergy} to identify
the openly $1/k$ geodesics.  First, we do not restrict ourselves
to $\Omega_k$ using an injectivity radius, nor do we restrict the
values of the energy.  This allows us to search for long and short
openly $1/k$ geodesics.  

\begin{defn} \label{unifen}
Let $E=E_{1/k,1/k,...1/k}:(M)^k \to \RR$ be called the uniform
energy.
\end{defn}

\begin{coro} \label{idenergycor}
For any openly $1/k$ geodesic $\gamma:S^1 \to M$ and any $t\in S^1$
the point 
\be
(\gamma(t), \gamma(t+2\pi/k), \gamma(t+4\pi/k),...\gamma(t-\pi/k))\in (M)^k
\ee
is a smooth critical point
of the uniform energy on $M^k$.
As we run through all values of $t$ we get a critical level set,
which is itself a closed geodesic in $(M)^k$.
\end{coro}

Before we set up the converse, we add a short lemma about geodesics
generated by critical points.

\begin{lemma} \label{idenergy1}
If $\bar{x}=(x_1,...x_k)\in (M)^k$ is a smooth critical point of 
the uniform energy $E:M^k \to \RR$,
then it defines a unique closed geodesic, $\gamma_{\bar{x}}$, 
which runs minimally between
the cyclic permutations $(x_1,x_2,...x_k)$, $(x_2, x_3,...x_k, x_1)$,
$(x_3,...x_k, x_1, x_2)$ and finally back through $(x_k, x_1,...x_{k-1})$
to $(x_1,...x_k)$.
\end{lemma}

\Pf
We know from Theorem~\ref{oldenergy} that if $\bar{x}$ is a critical
point we get a unique geodesic $\gamma:S^1\to M$ running through
$x_1, x_2,$ and on through $x_k$ and back to $x_1$.  So we can just
take 
\be
\bar{\gamma}(t)=(\gamma(t), \gamma(t+2\pi/k),...\gamma(t+(k-1)\pi/k)).
\ee
\qed

\begin{defn}
If $\bar{x}$ is a smooth critical point such that every point on
$\gamma_{\bar{x}}$ is also a smooth critical point, then we
say $\bar{x}$ is a rotating smooth critical point and 
$\gamma_{\bar{x}}$ is a rotating smooth critical set.
\end{defn}

\begin{theorem} \label{openenergy}
Openly $1/k$ geodesics in a convex compact Riemannian manifold with boundary,
$M$, have a one to one correspondance
with rotating smooth critical points of the uniform
energy in $(M)^k$ of nonzero value.
\end{theorem}

\Pf
This pretty much follows from Theorem~\ref{oldenergy},
Corollary~\ref{idenergycor} 
and Lemma~\ref{idenergy1}.
\qed

\begin{coro} \label{openenergyc} 
Given a manifold $M$, it's open index, $opind(M)$, is the
smallest value $k$ such that uniform energy $E: M^k \to \RR$ has a 
rotating smooth critical point with a nonzero value.
\end{coro}

\begin{example} \label{openminind0}
Suppose we use this approach to study the length
spectrum of $S^1$.  First we verify that
$L^{open}_{1/2}(S^1)=\emptyset$
because
\be
E(s,t)= 4(|s-t|\, mod \, 2\pi)^2= 4(s-t)^2 \, mod \, 16\pi^2
\ee
has only $(0,0)$ as a smooth critical point.
For $L^{open}_{1/3}(S^1)=\{2\pi\}$ we examine
\be
E(s,t,r) =3 (|s-t|\, mod \, 2\pi)^2+3(|t-r|\, mod \, 2\pi)^2+3(|r-s|\, mod \, 2\pi)^2
\ee
This energy is smooth as long as $|s-t|\neq 2k\pi$, $|t-r|\neq 2k\pi$ and
$|r-s|\neq 2k\pi$.  For
$(s_0,t_0,r_0)$ in this domain, there are values $k_1, k_2, k_3\in \ZZ$
such that for all $(s,t,r)$ near $(s_o,t_0,r_0)$ we have:
\begin{eqnarray}
(s-t)\, mod \, 2\pi & =& s-t+2k_1\pi \\
(t-r)\, mod \, 2\pi &=& t-r +2k_2\pi \\
(r-s) \, mod \, 2\pi & = & r-s+2k_3 \pi
\end{eqnarray}
so
\be
E(s,t,r) =3 (s-t+2k_1\pi)^2+3(t-r+k_2\pi)^2+3(r-s+k_3\pi)^2 
\ee
Thus we can differentiate and get:
\begin{eqnarray}
0=\partial E/\partial s &=& \, 2 (s-t +2 k_1 \pi) -2(r-s +2 k_2 \pi) \\
0=\partial E/ \partial t&=&  -2(s-t +2 k_2 \pi)+2(t-r +2 k_3 \pi) \\
0=\partial E/ \partial r&=& \, -2(t-r +2 k_3 \pi)+2(r-s +2 k_1 \pi) 
\end{eqnarray}
which implies that 
\be
h=(s-t) \, mod \, 2\pi = (t-r) \, mod \, 2\pi = (r-s) \, mod \, 2\pi.
\ee
Since $3h mod \, 2\pi =0$ we know
our smooth critical points have the form $(s,s,s)$ or
$(s, s+2\pi/3, s + 4\pi/3)$ or
$(s, s+4\pi/3, s + 2\pi/3)$.
This gives us two nonzero rotating critical points whose
energy is $9 (2\pi/3)^2=4\pi^2$, so their length is $2\pi$.
Thus $L^{open}_{1/3}=\{2\pi\}$.
Thus $openind(S^1)=3$.
\end{example}

Using a similar analysis of other compact length spaces one should
be able to impose lower bounds on their
minimizing index [Problem~\ref{openminind}].

Lusternick and Fet proved the existence of closed geodesics
on an arbitrary compact Riemannian manifold by producing
critical points of the energy functional.  Such critical
points are produced using Morse Theory and the topological properties
of the product space.  It is much more difficult to prove the
existance of rotating critical points.  [Problem~\ref{rotcritptsq}]
In fact, not all compact length spaces 
have closed geodesics.

\begin{example} \label{emptyL0}
Let $X=[0,1]$ with the standard metric $d(s,t)=|s-t|$.
Then for any $k \in \NN$, we study
\be
E(s_1, s_2,...s_k)=\sum_{j=1}^k k(s_j-s_{j+1})^2 \textrm{ where }s_{k+1}=s_1.
\ee
This is a smooth function on $(0,1)^k\subset [0,1]^k$, and its
critical points satisfy 
\be \label{int1}
s_{j-1}-s_j=s_j-s_{j+1} \textrm{ for } j=1,...k.
\ee
Since we are not on a circle, these points cannot wrap around, so
(\ref{int1}) implies that all the $s_j=0$.
Thus there are no smooth critical points and by Theorem~\ref{openenergy}
$X=[0,1]$ has no openly $1/k$ geodesics for any $k$ and 
by Theorem~\ref{openbigcup} it has
no closed geodesics at all.
\end{example}

In a similar manner Theorems~\ref{openenergy} and~\ref{openbigcup}
could be used to prove other compact length spaces have no
closed geodesics.  [Problem~\ref{emptyL}]

\begin{figure}[htbp]
\centering
\includegraphics[height=1in ]{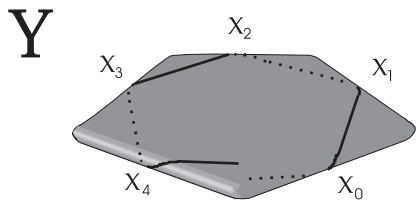}
\caption{The points $x_i \in Y$ correspond to a critical point 
$(x_0, x_1, x_2, x_3, x_4)$ of the uniform energy function on $(Y)^5$.
The geodesic $\gamma$ here is approaching $x_0$ from above.}
\label{pentagonfig}
\end{figure}

\begin{rmk} \label{nonsmooth}
Naturally, one would like to extend
Theorem~\ref{openenergy} to obtain some method of detecting
a $1/k$ geodesic which is not openly $1/k$.  To do so, one
might consider selecting nonsmooth critical points using
techniques from Grove-Shiohama's
critical points of distance function  or Chang's
critical points of Lipschitz functions \cite{GrShio} \cite{Cng}.

Using such techniques one would detect the $1/4$ geodesic
in Example~\ref{chiclet} (see Figure~\ref{chicletfig}.  That is
the point $(p_0, p_1, p_2, p_3)\subset (M_\epsilon)^4$ defined
using the $p_i$ in Example~\ref{chiclet} is such a nonsmooth
critical point.

Similarly, if one were to take a tubular neighborhood of a solid
pentagon, $Y$, 
in the $xy$ plane instead of a square as in Figure~\ref{pentagonfig}
and look at five evenly spaced points, $x_j$, on the equator 
near the midpoints of the sides of the pentagon, then 
one would again get a nonsmooth critical point in
the sense of Chang or of Grove-Shiohama.

However, if we let
$\gamma$ run minimally with positive $z$ from $x_0$ to $x_1$
and minimally from $x_1$ to $x_2$, one could verify it was running
minimally from $\gamma(t)$ to $\gamma(t+d(x_0,x_1))$, just like
the squarelike geodesic in Example~\ref{chiclet}.
However, if we continue to extend $\gamma$ in this manner
alternating above and below, it returns to $\gamma(0)$ from
above creating a corner!  So there is no geodesic
corresponding to this nonsmooth critical point, although
it is halfway around a $1/10$ geodesic. %ADD FIGURE

The author proposes in Problem~\ref{nonsmoothq} to study
nonsmooth critical points.
\end{rmk}

%\psdraft
%\psfull

\begin{rmk}\label{nondegen1}
An advantage of focusing on smooth critical points is that we can 
discuss the Hessian of the energy and degeneracy.
Naturally each openly $1/k$ geodesic is a degenerate critical
point because of the fact that there is an entire critical
level $\gamma_{\bar{x}}$.
However, a closed geodesic is said to be ``degenerate'' iff the 
$|det Hess _\perp E|=0$ where we focus on the directions perpendicular to
this rotational degeneracy.  Such geodesics have smoothly closed
Jacobi fields perpendicular to $\gamma'$. \cite{BTZ}\cite{GlMy}
In fact, there should be a stronger more global statement
describing an openly $1/k$ geodesic which corresponds to
a nondegenerate critical point of $E:M^k \to \RR$ [Problem~\ref{nondegenq}].
\end{rmk}

On smooth Riemannian manifolds the Morse index is the index of the Hessian 
of the energy of a geodesic, $\int |\gamma'(t)|^2\,dt$.
In particular, index of a closed geodesic, denoted $ind(\gamma)$, 
is the dimension
of the subspace of smooth closed vector fields perpendicular
to $\gamma'$, $V_\lambda$,
on which $H$ is negative definite, where
\be
H(X,Y)=\int_0^{2\pi}
<\nabla X, \nabla Y> -<R(X, \gamma'(t)) \gamma'(t), Y> \,dt.
\ee 
Morse proved that for geodesic segments, where the vector fields
have no assumption on periodity, the index bounds the number of
conjugate points on a segment.  Closed Geodesics have been studied
by Klingenberg and Ballman-Thorgbergson-Ziller, relating their index
to the Poincare Map \cite{K2} \cite{BTZ}.%\cite{GlMy} 

It is important to emphasize here that the Morse Index is
defined using vector fields and a covariant derivative and thus is not
naturally extended to compact length spaces.  Even when viewed as
a Hession of an energy on the loop space there is a significant
difficult defining an extension of the concept.  Finally, the Poincare
Map and even the unique extension of a geodesic is not defined
on arbitrary compact length space.  

\begin{theorem} \label{morsethm}
The Morse Index of a geodesic, $\gamma$, in a compact
Riemannian manifold satisfies: 
\be \label{morsethm1}
ind(\gamma) \le (n-1)(opind(\gamma)).
\ee
\be  \label{morsethm2}
ind(\gamma) \le (n-1)(minind(\gamma)+1).
\ee
\end{theorem}

\Pf
Since (\ref{morsethm1}) and Lemma~\ref{minop} imply 
(\ref{morsethm2}), we can concentrate on an openly
$1/k$ geodesic, $\gamma$.

Let $t_j=2\pi j/k$ for $j=0$ to $k-1$.
Following \cite{BTZ}, we have
$V_\Lambda$ equal to the direct sum of $V_\Lambda^1$ and $V_\Lambda^2$
where $V_\Lambda^1$ are piecewise Jacobi along this partition
and $V_\Lambda^2$ are smooth vector fields $=0$ on the partition.
They are orthogonal with respect to $H$ and $H$ is positive
definite on $V_\Lambda^2$ because the geodesic is minimal
between the points on the partition.  Note that the crucial point
is that we do not use the injectivity radius here.  Instead the number
of points in the partition depends on the openly
minimizing index of $\gamma$.
This immediately proves that the Morse Index of $\gamma$ satisfies:
\be 
ind(\gamma) \le dim V_\Lambda^1 = (n-1)k.
\ee
\qed

\begin{example}
A $1/k$ geodesic may have Morse index $0$ no matter how large $k$ is,
as can be seen in spaces with nonpositive sectional curvature,
like a torus, which have no conjugate points.
\end{example}

\begin{rmk}\label{hessunifen}
The crucial difference between the Morse Index and
the minimizing index of a closed geodesic is that the
Morse index is a purely local concept while the minimizing
index is a global concept checking for cut as well as
conjugate points.  

It would be interesting to investigate whether the
minimizing index of an openly $1/k$ geodesic is related
to the Hessian of the uniform energy in Theorem~\ref{openenergy}
[Problem~\ref{hessunifenq}].
\end{rmk}

\begin{rmk} \label{nondegen2}
Now if $M_i$ converge to $M$ in the $C^4$ sense their 
finite dimensional loop spaces $M_i^k$ converge in the $C^4$ sense.  
It is not hard to show (c.f. \cite{Cnly}) that a suddenly appearing critical
point under $C^4$ convergence must be degenerate.  Thus it is of
significant interest to identify these nondegenerate openly  $1/k$
geodesics.
\end{rmk}

\begin{rmk}\label{nondegen3}
One might be tempted to prove that the nondegenerate length spectrum
is a continuous function of smoothly converging manifolds.  However,
this can be seen not to be the case in Figure~\ref{boobah} since
although the nondegenerate length spectrum of the $M_i$ does converge
to the nondegenerate length spectrum of $Y$, there is a sequence of 
nondegenerate geodesics in the $N_i$ converging to the degenerate 
geodesic in $Y$.  In fact, Klingenberg proved that any geodesic
can be made into a limit of nondegenerate geodesics of a sequence of
$C^4$ close metrics on the manifold \cite{K2}.  Conley showed that
if the geodesic is degenerate then the sequence approaching it must
have a cancelling set of geodesics just as in Figure~\ref{boobah}
\cite{Cnly}.
\end{rmk}

%moved a section to almost.tex

%=======================

\sect{Open Problems} \label{sectproblems}

In this section we state some open problems, many of which
were mentioned earlier in the paper.  If you wish to work on 
one of these problems or have solved one, please
let the author know.

\begin{problem} \label{emptyL}
What compact length spaces have empty length spectra?
Luisternik and Fet proved that on any compact Riemannian manifold there
exist closed geodesics by proving the existence
of critical points of the energy functional
on the loop space (c.f. \cite{Cr}).
Here we need more than just critical points, so it would be
easier to prove some spaces have empty length spectra
using Theorem~\ref{openenergy}
in a manner similar to Example~\ref{emptyL0}.
\end{problem}

\begin{problem}
Are there upper bounds on $\min \, L_{1/k}(M)$ which depend on
volume rather than diameter?
See also Problem~\ref{volsysq}.
\end{problem}

\begin{problem} \label{minindMq1}
Find a compact manifold, $M$, whose shortest closed
geodesic has a larger minimizing index than the manifold.
%By Theorem~\ref{minindMinj}, $minind(M)$ is then 
%$\le \min L(M)/injrad(M)$.
\end{problem}

\begin{problem} \label{klingcat} 
Note that Klingenberg's Lemma implies
that the minimizing
index of any manifold without conjugate points is $2$  
[Corollary~\ref{klingcor}] . This includes all manifolds
with sectional curvature $\le 0$.  
What can one say about the minimizing index of a  CAT(0) spaces?
Suggestions for Problem~\ref{openminind} may help.
\end{problem}

\begin{problem}\label{conjpt}
Is there an appropriate definition for a conjugate point
on a compact length space which will give results
as strong as Corollary~\ref{klingcor}?
One might look at \cite{So1}, which has a definition
of conjugate point defined for an entirely different 
situation.  Keep Example~\ref{Klingitex} in mind.
\end{problem}

\begin{problem}\label{rotnabq}
Can one use the proofs of Rotman and Nabutovsky-Rotman's
results to provide bounds on $\min L_{1/k}(M)$?
See Proposition~\ref{Ro} nor~\ref{rotnab}.
\end{problem}

\begin{problem} \label{volsysq}
Try to extend the volume estimates on $\min L_{1/2}$
given in Propositions~\ref{BavSak}, ~\ref{Low1}, ~\ref{Pu}
and~\ref{Ro}
to compact length spaces with finite second Hausdorff measure.
It would not be expected that the results would follow without
some additional conditions.
See Remark~\ref{volsys}.
\end{problem}

\begin{problem}\label{whennonempty}
What properties can be placed on a simply connected manifold to guarantee the
existence of a $1/2$ geodesic?
Note that Theorem~\ref{openenergy} cannot be used to find
a $1/2$ geodesic but Problem~\ref{nonsmoothq} might prove
helpful.  
\end{problem}

\begin{problem} \label{boundminind}
What properties can be placed on a manifold to allow
one to estimate its minimizing index?
See Problem~\ref{openminind} for one possible approach.
\end{problem}

\begin{problem} \label{minopind}
In Lemma~\ref{minop} we related the open index to the
minimizing index of a geodesic.  On the standard sphere 
the difference between these indices is exactly 1
for all geodesics.  What other manifolds share this property?
[c.f. 
%Theorem~\ref{morseop}, 
% FIX Proposition~\ref{klap}
and Theorem~\ref{openenergy}].
\end{problem}

\begin{problem} \label{minindMdq1}
Is there a version of Proposition~\ref{minindMd} which involves
the volume rather than the diameter of the manifold?
\end{problem}

\begin{problem}  \label{minindMdq2}
What happens in the equality case for Proposition~\ref{minindMd}?
\end{problem}

\begin{problem} \label{exact}
On a manifold with minimizing index, $minind(M)=k$, is there
an exact bound on $\min \, L(M)$ which depends on $k$?  
\end{problem}

\begin{problem} \label{bangertq}
If $M_i$ converge to $S^2$ with the standard metric in the
Gromov-Hausdorff sense, and they have a common lower bound
on their injectivity radius, $injrad(M_i)\ge i_0>0$, then
do all prime geodesics $\gamma_i:S^1\to M_i$ satisfy
Bangert's Theorem that $L(\gamma_i)$ either converge
to $2\pi$ or diverge to infinity?  See Remark~\ref{bangertrmk}.
\end{problem}

\begin{problem}\label{rotcritptsq}
Many theorems proving the existence of a closed geodesic
on a manifold involve the study of the Morse Theory of
the loop space and the existence of critical point on
that loop space.  
To produce a $1/k$ geodesic, Theorem~\ref{openenergy}
requires that we find a ``rotating'' critical point
of an energy on a product space.   What conditions can be
placed on a manifold or metric space to prove the existence
of such a critical point?
\end{problem}

\begin{problem}\label{openminind}
Estimate the minimizing index of a compact length
space or provide a lower bound on the minimizing index
using Theorem~\ref{openenergy}.  See Example~\ref{openminind0}
for a simple case.  Such an estimate would then provide
an  estimate on $minind(M)$ and
thus $\min L(M)$ via Theorem~\ref{minindMd} and Lemma~\ref{minop}.
\end{problem}

\begin{problem}\label{nonsmoothq}
In Theorem~\ref{openenergy}, we relate openly $1/k$ geodesics
on a compact Riemannian manifold $M$
to special smooth critical points of an energy on
$M^k=M\times M\times \cdots \times M$.  It would be interesting
to study whether some definition for a nonsmooth
critical point might be used that relates to
$1/k$ geodesics. See Remark~\ref{nonsmooth}.
This might help solve Problem~\ref{whennonempty}.
\end{problem}

\begin{problem} \label{nondegenq}
A degenerate closed geodesic is a geodesic whose
energy functional is degenerate.  It has been proven
to have a smoothly closed Jacobi field in \cite{BTZ}\cite{GlMy}.
Is there a similar more global property concerning
nearby geodesics for a degenerate openly $1/k$ geodesic
where one defines degenerate using the Hessian of
the uniform energy?
See Remarks~\ref{nondegen1}, ~\ref{nondegen2} and~\ref{nondegen3}.
\end{problem}

\begin{problem}\label{hessunifenq}
Does the index of the Hessian of the uniform energy
provide an estimate on the minimizing index?
It would be interesting to investigate whether the
minimizing index of an openly $1/k$ geodesic is related
to the Hessian of the uniform energy in Theorem~\ref{openenergy}.
See Remark~\ref{hessunifen}.
\end{problem}

\begin{problem} \label{volgapq} 
In Theorem~\ref{volgap}, we estimate the location of the
length spectrum of a Riemannian manifold, $N^n$, whose volume
is close to that of the sphere and whose Ricci curvature
is bounded from below.  Can one find an explicit formula
for the estimating function, $\Psi$? How strong is
its dependance on $k$?   Can one control $L(N^n)$ and
not just $L_{1/k}$?  Note that Colding's Volume Theorem \cite{Co1}
does not give a precise estimate on the Gromov-Hausdorff
convergence and getting one from his proof would be very difficult.
However, proving this result directly may be possible.
See Remark~\ref{volgapr}.
\end{problem}

\begin{problem}\label{volgapmin}
Find a sequence of manifolds $M_j^n$ with $Ricci \ge (n-1)$,
$Vol(M^n_j) \to Vol(S^n)$ such that $\min L(M^n_j) \to 0$
or prove this cannot occur.  Note in Example~\ref{football}
we showed there is no uniform lower bound on injectivity radius
implies by the $Ricci$ and volume conditions.
\end{problem}

\begin{problem}\label{volgapdis}
Find a sequence of manifolds $M_j^n$ with $Ricci \ge (n-1)$,
$Vol(M^n_j) \to Vol(S^n)$ with $L_j \in L(M_j)$
such that $L_j\to L_\infty \notin L(S^2)$ or prove this
cannot occur.  Note that by Theorem~\ref{volgap}
we know the $\gamma_j$ of length $L_j$ have $minind(\gamma_j)\to\infty$.
It is quite possible that the $M_j$ in Example~\ref{football}
have disappearing geodesics, so these surfaces are worth investigation.
One might begin by stretching elastic loops around footballs
in a clever way.
\end{problem}

\begin{problem} \label{diamgapq}
Is it possible to get a gap theorem for manifolds with
$Ricci\ge (n-1)$ and diameter close to $\pi$?  
See Remark~\ref{diamgap}.
\end{problem}

\begin{problem}\label{sphsusp}
Given a length space $X$ what can one say about the length
spectrum of the spherical suspension over $X$?
See \cite{BBI} for a rigorous definition of a spherical
suspension.
\end{problem}

\begin{problem}\label{radgapmin}
Find a sequence of manifolds $M_j^n$ with $Ricci \ge (n-1)$,
$rad(M^n_j) \to rad(S^n)$ such that $\min L(M^n_j) \to 0$
or prove this cannot occur.  See Remark~\ref{footballrad}.
\end{problem}

\begin{problem}\label{radgapdis}
Find a sequence of manifolds $M_j^n$ with $Ricci \ge (n-1)$,
$rad(M^n_j) \to rad(S^n)$ with $L_j \in L(M_j)$
such that $L_j\to L_\infty \notin L(S^2)$ or prove this
cannot occur.  Note that by Theorem~\ref{radgap}
we know the $\gamma_j$ of length $L_j$ have $minind(\gamma_j)\to\infty$.
\end{problem}

\begin{problem}\label{torusgapmin}
Find a sequence of manifolds $M_j^n$ with $Ricci \ge -\epsilon_j(n-1)\to 0$,
and $b_1(M^n)=n$ such that 
$\min L(M^n_j) \to 0$
or prove this cannot occur.  See Theorem~\ref{torusgap}.
\end{problem}

\begin{problem}\label{torusgapdis}].
Find a sequence of manifolds $M_j^n$ with $Ricci \ge -\epsilon_j(n-1)\to 0$,
and $b_1(M^n)=n$ such that 
such that $L_j\to L_\infty \notin L(S^2)$ or prove this
cannot occur.  Note that by Theorem~\ref{torusgap}
we know the $\gamma_j$ of length $L_j$ have $minind(\gamma_j)\to\infty$.
\end{problem}

\begin{problem}\label{cosmosq}
Analyze the length spectra of locally almost isotopic manifolds
mentioned in Remark~\ref{cosmos}.
\end{problem}

\begin{problem} \label{PLBS}
In Section~\ref{almrig} we explained how some rigidity theorems
with extremal diameters, volumes or eigenvalues relative
to Ricci curvature bounds have stability statements.
Propositions~\ref{Pu}, \ref{Low1} and~\ref{BavSak} 
do not involve Ricci curvature but do have rigidity
results when their equalities have been achieved.
Do they have related stability or stability theorems?
Without the Ricci curvature bounds one wouldn't expect
these theorems to involve Gromov-Hausdorff convergence.
\end{problem}

\end{document}